\documentclass[10pt]{amsart}
\usepackage{amsmath}
\usepackage{amssymb}
\usepackage{amsthm}
\usepackage{amsfonts}
\usepackage{amstext}
\usepackage{amsopn}
\usepackage{amsxtra}
\usepackage{pst-all}
\usepackage[latin1]{inputenc}
\usepackage{dsfont}
\newcommand\R{{\ensuremath {\mathbb R} }}
\newcommand\C{{\ensuremath {\mathbb C} }}
\newcommand\N{{\ensuremath {\mathbb N} }}
\newcommand\Z{{\ensuremath {\mathbb Z} }}

\newcommand\1{{\ensuremath {\mathds 1} }}
\renewcommand\phi{\varphi}

\newcommand{\gH}{\mathfrak{H}}

\newcommand{\gS}{\mathfrak{S}}
\newcommand{\wto}{\rightharpoonup}
\renewcommand{\to}{\rightarrow}
\newcommand{\cP}{\mathcal{P}}
\newcommand{\cN}{\mathcal{N}}

\newcommand{\cB}{\mathcal{B}}

\newcommand{\cF}{\mathcal F}

\newcommand{\cM}{\mathcal M}

\newcommand{\tr}{{\rm tr}\,}
\newcommand{\cE}{\mathcal{E}}

\newcommand{\cQ}{\mathcal{Q}}
\newcommand{\cK}{\mathcal{K}}

\newcommand{\cG}{\mathcal{G}}
\newcommand{\cX}{\mathcal{X}}
\newcommand{\Ex}{{\rm Ex}}
\newcommand\ii{{\ensuremath {\infty}}}
\newcommand\pscal[1]{{\ensuremath{\left\langle #1 \right\rangle}}}
\newcommand{\norm}[1]{ \left| \! \left| #1 \right| \! \right| }

\newcommand{\Hhalf}{H^{1/2}(\R^3)}

\newcommand{\IHF}{I^{\rm HF}}
\newtheorem{thm}{Theorem}
\newtheorem{lemma}{Lemma}
\newtheorem{corollary}{Corollary}
\newtheorem{prop}{Proposition}

\newtheorem{remark}{Remark}

\numberwithin{equation}{section}
 \numberwithin{lemma}{section}
\numberwithin{prop}{section}
\numberwithin{corollary}{section}
\renewcommand{\tr}{{\rm Tr} }
\newcommand{\str}{{\tr_{P^0_-}} }

\long\def\symbolfootnote[#1]#2{\begingroup%
\def\thefootnote{\fnsymbol{footnote}}\footnote[#1]{#2}\endgroup}

\renewcommand{\leq}{\leqslant}
\renewcommand{\geq}{\geqslant}

\begin{document}

\title[Minimizers for HFB Theory of Neutron Stars]{Minimizers for the Hartree-Fock-Bogoliubov \\ Theory of Neutron Stars and White Dwarfs}

\author[E. Lenzmann]{Enno LENZMANN}
 \address{Institute for Mathematical Sciences, University of Copenhagen, Universitetsparken 5, 2100 Copenhagen \O, Denmark.}
 \address{Department of Mathematics, M.I.T., 77 Massachusetts Avenue, 02139 Cambridge, MA. USA.}  
  \email{lenzmann@math.ku.dk}

\author[M. Lewin]{Mathieu LEWIN}
 \address{CNRS \& Laboratoire de Mathématiques (CNRS UMR 8088), Université de Cergy-Pontoise, 95000 Cergy-Pontoise, France.}
  \email{Mathieu.Lewin@math.cnrs.fr}

\date{March 20, 2010. Final version to appear in {\em Duke Math.~Journal}}

\begin{abstract}
We prove the existence of minimizers for Hartree-Fock-Bogoliubov (HFB) energy functionals with attractive two-body interactions given by Newtonian gravity. This class of HFB functionals serves as model problem for self-gravitating relativistic Fermi systems, which are found in neutron stars and white dwarfs. Furthermore, we  derive some fundamental properties of HFB minimizers such as a decay estimate for the minimizing density. 

A decisive feature of the HFB model in gravitational physics is its failure of weak lower semicontinuity. This fact essentially complicates the analysis compared to the well-studied Hartree-Fock theories in atomic physics. 
\end{abstract}

\maketitle

\setcounter{tocdepth}{1}
\tableofcontents

\section{Introduction}

The Hartree-Fock-Bogoliubov (HFB) theory is a widely used tool \cite{RinSch-80,DeaHjo-03}  for understanding many-body quantum systems where attractive two-body interactions are dominant. In particular, HFB energy functionals incorporate the physical phenomenon of {\em ``Cooper pairing''} which is likely to occur whenever attractive forces among quantum particles become significant; e.\,g., in nuclear and gravitational physics as well as in superconducting materials. Despite the broad range of physical applications for HFB theory, not much has been known rigorously concerning the existence of minimizers, let alone a proof of their fundamental properties such as Cooper pair formation. 

As a starting point for rigorous analysis, this paper is devoted to the existence of minimizers for HFB relativistic energy functionals with an interaction which behaves at infinity like the (attractive) Newtonian gravity. This specific class of functionals can be viewed as a model problem for self-gravitating relativistic Fermi systems which are found, for example, in {\em neutron stars} and {\em white dwarfs.} Moreover, speaking from a mathematical point of view, the corresponding variational problem exhibits the delicate property of criticality, as we will detail below.

The most challenging main feature of the HFB variational problem in gravitational physics is its lack of weak lower semicontinuity (wlsc) due to the attractive interaction among particles. As a consequence of the absence of wlsc, the existence proof for minimizers is much more involved than for the well-studied Hartree-Fock (HF) models arising in atomic physics, where wlsc plays an essential role; see \cite{LieSim-77,Lions-87}.  Another difficulty in the analysis of HFB models stems from its translational invariance and the fact that the main variables are {\em two operators} related via a complicated constraint inequality. Finally, a further complication (although conceptually less important) is the treatment of the pseudo-differential operator describing the kinetic energy of relativistic fermions.

\medskip

Following the general rigorous discussion of HFB energy functionals in \cite{BacLieSol-94}, we specifically consider the energy functional given by
\begin{multline}
\cE(\gamma,\alpha) :=\tr \left ( T \gamma \right ) + \frac{\kappa}{2} \iint_{\R^3\times\R^3}W(x-y)\rho_\gamma(x)\rho_\gamma(y)dx\,dy\\ - \frac{\kappa}{2} \iint_{\R^3\times\R^3}W(x-y)|\gamma(x,y)|^2dx\,dy
+\frac{\kappa}{2} \iint_{\R^3\times\R^3}W(x-y)|\alpha(x,y)|^2dx\,dy ,
\label{def_energy_intro}
\end{multline}
which models a system of relativistic fermions subject to the interaction $W$. We will assume that $W(x)$ behaves like $-1/|x|$ at infinity, i.e. that far away the attractive Newtonian force is predominant\footnote{See Remark \ref{rmk:general_potential} for precise assumptions on $W$.}. Indeed in most of the paper we will even take for simplicity 
$$W(x)=-\frac1{|x|}.$$
In \eqref{def_energy_intro}, $\tr (\cdot)$ denotes the trace, and the pseudo-differential operator  
$$T=\sqrt{-\Delta + m^2}-m$$
describes the kinetic energy of a fermion with mass $m > 0$. The coupling constant $\kappa >0$ parametrizes the strength of the interaction among the particles. Furthermore, the variables $\gamma$ and $\alpha$ are two \emph{operators} acting on $L^2(\R^3; \C^2)$.\footnote{More generally, we shall consider $L^2(\R^3; \C^q)$ below, where $q \geq 1$ denotes the internal spin degree of freedom. The case $q=2$ (corresponding to spin 1/2) is physically the most relevant one.} The operator $\gamma$ is referred to as the \emph{one-body density matrix}. It is self-adjoint, nonnegative and has a finite trace fixed to be
\begin{equation}
\tr(\gamma)=\lambda,
\label{constraint_1} 
\end{equation}
for some given $\lambda >0$, which is a real parameter and it corresponds to the expected value of the number of particles in the star. The density $\rho_\gamma$ is the unique $L^1(\R^3)$ non-negative function that satisfies $\tr(\gamma V)=\int_{\R^3}\rho_\gamma(x)V(x)\,dx$ for any bounded function $V$. As $\gamma$ is trace-class, it has a kernel $\gamma(x,y)$ appearing in the second line of \eqref{def_energy_intro}, which is a $2\times 2$ hermitian matrix for almost every $x,y\in\R^3$. Finally, the operator $\alpha$ is called the \emph{pairing density matrix}. It is only assumed to be Hilbert-Schmidt, that is we have $\tr(\alpha^*\alpha)<\ii$. Its kernel is a $2\times 2$ matrix which is supposed to be antisymmetric in the following sense: $\alpha(x,y)^T=-\alpha(y,x)$, where $^T$ is the usual transposition of matrices.

\medskip
To formulate the HFB variational problem for the energy $\cE(\gamma,\alpha)$, we have to supplement the side condition (\ref{constraint_1}) by the following operator inequality relating $\gamma$ and $\alpha$:
\begin{equation}
\left(\begin{matrix}
0 & 0\\
0 &  0
\end{matrix}\right)\leq \left(\begin{matrix}
\gamma & \alpha\\
\alpha^* &  1-\overline{\gamma}
\end{matrix}\right)\leq \left(\begin{matrix}
1 & 0\\
0 &  1
\end{matrix}\right)\quad\text{on } L^2(\R^3;\C^2)\oplus L^2(\R^3;\C^2).
\label{constraint_2}
\end{equation}
This inequality guarantees that the pair $(\gamma,\alpha)$ is associated to a unique quasi-free state in Fock space \cite{BacLieSol-94}. The corresponding HFB minimization problem then reads
\begin{equation*}
 \boxed{ I(\lambda):=\inf\big\{\cE(\gamma,\alpha) : \gamma^*=\gamma,\ \alpha^T=-\alpha \ \text{such that \eqref{constraint_1} and \eqref{constraint_2}}\big\}.}
\end{equation*}
For the physical interpretation and some background from many-body quantum mechanics, we refer the reader to \cite{BacLieSol-94,BlaRip-85}.

\medskip
Our first main result is formulated in Theorem \ref{thm:relaxed} below, which shows that
$$ \boxed{ \mbox{$I(\lambda)$ is attained for $0 \leq \lambda < \lambda^{\rm HFB}(\kappa)$ and $\kappa < 4/\pi$.} } $$
More precisely, we prove that all minimizing sequences (with such $\lambda$ and $\kappa$) are relatively compact, up to translations, in an appropriate topology for the pair of operators $(\gamma, \alpha)$. In physical terms, the finite number $\lambda^{\rm HFB}(\kappa)$ can be interpreted as the {\em Chandrasekhar limiting mass} for the HFB model considered here. That is, a self-gravitating relativistic Fermi system with supercritical particle number, 
$$
\lambda > \lambda^{\rm HFB}(\kappa),
$$ 
can undergo {\em ``variational collapse''} in the sense that  $\cE(\gamma,\alpha)$ fails to be bounded below and hence minimizers cannot exist in this case. Regarding the relation to semi-classical models of stars, we remark that a useful estimate for $\lambda^{\rm HFB}(\kappa)$ is provided by the asymptotic estimate (when the interaction is purely gravitational):
\[
\lambda^{\rm HFB}(\kappa) \sim \left (\frac{\tau_c}{2^{1/3}\kappa} \right )^{3/2} \quad \mbox{as} \quad \kappa \to  0. 
\]
Here $\tau_c=O(1)$ is a numerical constant that can be calculated from the classical {\em Thomas-Fermi-type theory} for neutron stars and white dwarfs, which was initiated by S.~Chandrasekhar in his seminal work \cite{Chandra-31}; see also \cite{LieThi-84,LieYau-87}.

\medskip
Apart from the existence of minimizers, we also establish some of their fundamental properties. In particular, we derive a certain decay estimate for the minimizing density $\rho_\gamma(x)$, which turns out to be essential when proving the existence of minimizers, as we will detail below. Here, we just point out that the fact that a decay estimate for $\rho_\gamma$ {\em cannot} be simply obtained by using the self-consistent equation for HFB minimizers, since the minimizing one-particle matrix $\gamma$ is not known to have a finite rank. In fact, by using results from \cite{BacFroJon-08,FraLieSeiSie-07}, we shall see below that any minimizer $(\gamma,\alpha)$ in the spin-1/2 case exhibits an infinite rank, provided we assume $\alpha \neq 0$ holds. To deal with this difficulty arising from infinite-rank properties of HFB minimizers, we devise some variational arguments combined with a second order expansion of the HFB energy to prove an estimate for the fall-off behavior of the minimizing density $\rho_\gamma(x)$.

Concerning further fundamental properties of minimizers for $\cE(\gamma, \alpha)$, we remark that it is an important open question to show that minimizers actually exhibit non-vanishing pairing $\alpha \neq 0$, at least for a coupling constant $\kappa$ which is not too small (relatively to $\lambda^{-2/3}$). On heuristic grounds, one expect such a phenomenon of {\em ``Cooper pair formation''} to be energetically favorable due to the attractive interaction among particles. However, it seems to be a formidable task to find mathematical proof for this claim.  

\medskip

Let us briefly comment on our existence proof. The relative compactness of minimizing sequences (and, in particular, the existence of minimizers) is the consequence of the validity of the \emph{binding inequality:}
\begin{equation}
 I(\lambda)<I(\lambda')+I(\lambda-\lambda'), \quad \mbox{for all $\lambda' \in (0, \lambda)$,}
\label{binding_gHF}
\end{equation}
which is also referred to as {\em strict sub-additvity condition}. In fact, binding inequalities like \eqref{binding_gHF} appear naturally in the analysis of compactness properties of minimizing sequences; for instance, when using the ``concentration-compactness principle'' as done in \cite{Lions-84,Lions-84b}. Moreover, binding inequalities also turn out to be useful for linear models, in which the bottom of the essential spectrum has the form of the minimum on the right side in \eqref{binding_gHF}, as expressed by the celebrated HVZ Theorem; see \cite{Hun-66,VanWinter-64,Zhislin-60,Enss-77,Simon-77,Sigal-82}. The interpretation of \eqref{binding_gHF} is that it is not favorable for a minimizing sequence to have a mass $\lambda-\lambda'$ escaping to infinity, while the mass $\lambda'$ stays in a bounded set (up to translations).

However, there is a notable difficulty if one tries to prove \eqref{binding_gHF} for the HFB model studied in this paper. To see this, we recall that in the usual $N$-body Schr\"odinger setting, the binding inequalities are quantized, i.\,e., they are of the form
$$I_{\rm Sch}(N)< I_{\rm Sch}(N-N')+I_{\rm Sch}(N'),\quad \mbox{for all $N'=1\ldots N-1$} .$$
This quantization just reflects the fact that we can restrict ourselves to states with a definite particle number and use adapted geometrical methods \cite{Enss-77,Simon-77,Sigal-82}. In particular, an important feature of quantized binding inequalities is that they can be proved by induction; see, e.\,g., \cite{Friesecke-03}. In HFB theory, by contrast, it is not expected that ground states will have a definite particle number. Therefore, the binding inequality \eqref{binding_gHF} cannot be simply established by induction in this model.

To overcome this difficulty and eventually prove the binding inequality (\ref{binding_gHF}) for the HFB model considered here, we use a different method in the spirit of a technique introduced by P.-L.~Lions in \cite{Lions-87}. More precisely, we assume the existence of a non-relatively compact minimizing sequence with particle number $\lambda < \lambda^{\rm HFB}(\kappa)$, and we describe its behavior in detail. In particular, we show its splitting into several pieces receding from each other, where at least two parts are relatively compact and carry strictly positive particle numbers $\lambda^1 > 0$ and $\lambda^2 >0$, respectively. Furthermore, we conclude that the infima $I(\lambda^i)$ with $i=1,2$ must be attained, and we find that the ground state energy decomposes as
\begin{equation*}
I(\lambda) = I(\lambda^1) + I(\lambda^2) + I(\lambda- \lambda^1 - \lambda^2) .
\label{binding_eq}
\end{equation*}
On the other hand, by deriving a suitable decay estimate for any minimizer of $I(\lambda^i)$, we can deduce the strict inequality
\begin{equation*}
I(\lambda)< I(\lambda^1) + I(\lambda^2)  + I(\lambda- \lambda^1 - \lambda^2) ,
\label{binding_gHF2}
\end{equation*}
in contradiction to the previous equality. Thus the main assumption (failure of relative compactness of all minimizing sequences) cannot hold. This rationale establishes {\em a posteriori} the binding inequality (\ref{binding_gHF}) for the HFB model.

Although the aforementioned behavior of minimizing sequences appears physically quite clear, implementing the above argument mathematically is far from being an easy task. One difficulty stems from the fact that the main variable is a pair of \emph{two operators} satisfying the complicated constraint \eqref{constraint_2}, whereas in usual nonlinear PDE problems, the variable is often a single scalar or vector-valued function. In order to deal with this issue, we use some ideas of previous works \cite{HaiLewSer-05a,HaiLewSer-05b,FraLieSeiSie-07}. Also, we point out that the pairing term $\alpha$ cannot be handled by obvious means: In contrast to the trace-class operator $\gamma$, the operator $\alpha$ is \emph{a priori} only Hilbert-Schmidt. Moreover, loosely speaking, the direct and exchange term are both \emph{``subcritical''} in the sense that they may be controlled by the $L^{12/5}$-norm of $\sqrt{\rho}_\gamma\in H^{1/2}(\R^3)$. By contrast, the pairing energy depending on $\alpha$ is \emph{``critical''} because it can only be controlled by the kinetic energy of $\gamma$ itself.  A particular illustration of this difficulty is that even ruling out the vanishing of a minimizing sequence is quite delicate, and that the proof of the decay of minimizers is based on two-body techniques for the wavefunction $\alpha(x,y)$.

\subsection*{Organization of the Paper}
The outline of this paper is as follows. Section \ref{sec:basic} sets the stage for the variational calculus related to the HFB energy $(\gamma, \alpha)\mapsto\cE(\gamma, \alpha)$. Our main results are then formulated in Section \ref{sec:main}. After collecting some preliminary facts in Section \ref{sec:prelim}, the proofs of the main theorems are presented in Sections \ref{sec:proof_prop_min}--\ref{sec:proof_thm_orbitals} below. The appendix contains various technical results and proofs.

\subsection*{Notation} 
We employ standard notation for $L^p$ and Sobolev spaces. For an integer $q \geq 1$, we define the inner product $\langle f,g \rangle = \int_{\R^3} \sum_{i=1}^q \overline{f}_i(x) g_i(x) \,dx$ for functions $f$ and $g$ in $L^2(\R^3; \C^q)$. If $M$ is a $q\times q$ matrix, we use the notation $|M|^2=\tr(M^*M)$. In the whole paper we employ the physicists' notation $|\psi\rangle\langle\phi|$ for the operator $f\mapsto \pscal{\phi,f}\psi$. Furthermore, we shall use $\gS_p$, with $1\leq p<\ii$, to denote the Schatten class of operators $A$ acting on $L^2(\R^3;\C^q)$ and having a finite $p$-trace, and we endow it the norm $\norm{A}_{\gS_p}:=\tr(|A|^p)^{1/p}<\ii$; see \cite{Simon-79}. The space $\gS_\ii$ denotes the space of compact operators on $L^2(\R^3; \C^q)$, equipped with the operator norm $\| \cdot \|$. Below we shall also introduce a ``Sobolev-type'' Schatten space $\cX \subset \gS_1 \times \gS_2$ in the study of the HFB energy functional $\cE$.

For later convenience, we introduce two smooth cutoff functions $\chi_R$ and $\zeta_R$ that are defined as follows. Let $0 \leq \chi \leq 1$ be a fixed smooth function on $\R^3$ such that $\chi \equiv 1$ for $|x| < 1$ and $\chi \equiv 0$ for $|x| \geq 2$. For any $R > 0$, we then define the functions
\begin{equation} \label{eq:cutoff}
\chi_R(x) = \chi(x/R) \quad \mbox{and} \quad \zeta_R(x) = \sqrt{1-\chi_R(x)^2}. 
\end{equation}
Throughout this paper, we use the cutoff functions $\chi_R$ and $\zeta_R$ freely. Furthermore, we denote by $\tau_y : L^2(\R^3; \C^q) \to L^2(\R^3;\C^q)$ the unitary operator that is defined by 
\begin{equation} \label{eq:tau}
(\tau_y f) = f(\cdot-y) ,
\end{equation}
for given $y \in \R^3$. We also use the notation $U(\gamma,\alpha)U^* :=(U \gamma U^*,U \alpha U^*)$ for any unitary operator $U$ acting on $L^2(\R^3;\C^q)$, e.\,g., for the translation operator $U=\tau_y$. 

For the physically inclined reader, we remark that we work in units such that Planck's constant $\hbar$ and the speed of light $c$ satisfy $\hbar = c =1$.

\subsection*{Acknowledgments}
We are indebted to \'E.~S\'er\'e for helpful discussions when revising this paper. Also, we are grateful to P.~T.~Nam for his careful reading of the manuscript and for finding a gap in a previous version of this paper. We cordially thank V.~Bach, J.~Fr\"ohlich and L.~Jonsson for communicating their independent results \cite{BacFroJon-08} on the HFB model. 

The first author (E.\,L.) is partly supported by NSF Grant DMS-0702492 and a Steno research fellowship from the Danish science council. The second author (M.\,L.) is partly supported by ANR project ACCQuaRel. Finally, the authors gratefully acknowledge the kind hospitality of the Erwin Schr\"odinger Institut in Vienna, where parts of this work were completed.

\section{Basic Properties of HFB Energy} \label{sec:basic}

To prepare the statement of our main results, we first provide the adequate setting for the Hartree-Fock-Bogoliubov (HFB) variational problem, and we collect some basic properties needed for the rest of this paper. 

\subsection{HFB States With Finite Kinetic Energy}
To set up the variational calculus, we define a class of HFB states having finite pseudo-relativistic kinetic energy. Therefore we introduce the following (real) Banach space of density matrices
\begin{equation}
\cX = \big  \{(\gamma,\alpha) \in \gS_1\times\gS_2 \ : \ \gamma^* = \gamma,\ \alpha^T=-\alpha,\ \|(\gamma,\alpha) \|_{\cX}<\ii\big\},
\label{def_set_X}
\end{equation}
equipped with the norm
\begin{equation}
\|(\gamma,\alpha) \|_{\cX} = \big \| (1-\Delta)^{1/4} \gamma (1-\Delta)^{1/4} \big\|_{\gS_1} + \big \| (1-\Delta)^{1/4} \alpha \big \|_{\gS_2} .
\end{equation}
We remind the reader that $\gS_1$ and $\gS_2$ denote the space of trace-class and Hilbert-Schmidt operators on $L^2(\R^3; \C^q)$, respectively. The universally fixed integer $q \geq 1$ takes into account the internal spin degrees of freedom of the model. In what follows, we shall omit the dependence on $q$ whenever it is of no importance.

Furthermore, we define the following subsets of density matrices in $\cX$:
\begin{equation}
\cK =  \left\{(\gamma,\alpha) \in \cX\ :\ \left(\begin{matrix}
0&0\\ 0&0\end{matrix}\right) \leq \left(\begin{matrix}
\gamma&\alpha\\ \alpha^*&1-\overline{\gamma}\end{matrix}\right)\leq\left(\begin{matrix}
1&0\\ 0&1\end{matrix}\right) \right\}.
\end{equation}
and, for $\lambda > 0$ given, we set
\begin{equation}
\cK_\lambda = \left\{(\gamma,\alpha) \in \cX\ :\ \left(\begin{matrix}
0&0\\ 0&0\end{matrix}\right) \leq \left(\begin{matrix}
\gamma&\alpha\\ \alpha^*&1-\overline{\gamma}\end{matrix}\right)\leq\left(\begin{matrix}
1&0\\ 0&1\end{matrix}\right),\; \tr (\gamma) = \lambda \right\} .
\end{equation}
Here $\overline{\gamma}$ is defined by complex conjugation of the kernel $\gamma(x,y)$. Note that both $\cK$ and $\cK_\lambda$ are closed and convex sets in $\cX$, and that we have $0\leq\gamma\leq1$ for all $(\gamma,\alpha)\in\cK$. In the following, we shall use the \emph{notation:} 
\begin{equation}
\tr(T\gamma):=\tr\left((-\Delta+m^2)^{1/4}\gamma(-\Delta+m^2)^{1/4}\right)-m\tr(\gamma),
\end{equation}
for $T = \sqrt{-\Delta + m^2} -m$. Also, we sometimes allow the case of vanishing mass $m=0$. But we remind the reader that strict positivity of $m >0$ is an essential assumption in all the main theorems stated below.

\medskip
As a first simple fact, we obtain the following Sobolev estimate for the square root of the density function $\rho_\gamma$.
 
\begin{lemma} \label{lem:estim_density} For any $(\gamma,\alpha) \in \cK$, the density function $\rho_\gamma$ satisfies
\begin{equation*}
\big \langle \sqrt{\rho}_\gamma, \sqrt{-\Delta} \, \sqrt{\rho}_\gamma \big \rangle \leq \tr   (\sqrt{-\Delta} \, \gamma ) .
\end{equation*}
Moreover, the map $(\gamma, \alpha) \mapsto \sqrt{\rho}_\gamma$ is continuous from $\cK$ to $L^p(\R^3)$ for all $1\leq p<3$.
\end{lemma}

\begin{proof}[Proof of Lemma \ref{lem:estim_density}]
Writing $\gamma = \sum_n \lambda_n |\phi_n\rangle\langle\phi_n|$ with $0\leq\lambda_n\leq1$, we find
\begin{equation}
\tr \big (\sqrt{-\Delta} \, \gamma \big ) = \sum_n \lambda_n \langle \phi_n, \sqrt{-\Delta} \, \phi_n \rangle.
\end{equation} 
Furthermore, we recall that, for any $\phi, \psi \in \Hhalf$,
\begin{equation}
\big \langle \sqrt{|\phi|^2 + |\psi|^2}, \sqrt{-\Delta} \, \sqrt{|\phi|^2 + |\psi|^2} \big \rangle \leq \langle \phi, \sqrt{-\Delta}  \, \phi \rangle + \langle \psi, \sqrt{-\Delta} \, \psi \rangle .
\end{equation}
This follows from \cite[Theorem 7.13]{LieLos-01} and $\langle |\phi|, \sqrt{-\Delta} \, |\phi| \rangle \leq \langle \phi, \sqrt{-\Delta} \, \phi \rangle$. These two estimates are easily seen to imply the inequality in Lemma \ref{lem:estim_density}. The rest simply follows from Sobolev's embedding.\end{proof}

\subsection{Boundedness From Below of HFB Energy}
We define the (purely gravitational) HFB energy as
\begin{equation} 
\cE(\gamma,\alpha):=\tr \left ( T \gamma \right ) - \frac{\kappa}{2} D(\rho_{\gamma}, \rho_{\gamma} ) + \frac{\kappa}{2} \Ex(\gamma)
-\frac{\kappa}{2} \iint_{\R^3\times\R^3}\frac{|\alpha(x,y)|^2}{|x-y|}dx\,dy ,
\label{def_energy}
\end{equation}
where we now use the shorthand notations
\begin{equation}
D(f,g):=\iint_{\R^3\times\R^3}\frac{f(x)g(y)}{|x-y|}dx\,dy,\qquad \Ex(\gamma):=\iint_{\R^3\times\R^3}\frac{|\gamma(x,y)|^2}{|x-y|}dx\,dy,
\end{equation}
which we refer to as the {\em direct term} and {\em exchange term}, respectively. The last term on the right side of \eqref{def_energy} is called the {\em pairing term}.

For $(\gamma, \alpha) \in \cK$, we deduce from Lemma \ref{lem:estim_density} and the Hardy-Littlewood-Sobolev inequality that
\begin{equation}
D(\rho_\gamma,\rho_\gamma)\leq C \| \sqrt{\rho}_\gamma  \|^{4}_{L^{12/5}} < \infty .
\end{equation}
Thus the direct term $D(\rho_\gamma,\rho_\gamma)$ is well-defined. Next, we recall the pointwise estimate
\begin{equation}
|\gamma(x,y)|^2\leq\rho_\gamma(x)\rho_\gamma(y), \quad \mbox{for a.\,e.~$(x,y) \in \R^3 \times \R^3$},
\label{estim_exchange}
\end{equation}
which is obtained by writing the spectral decomposition of $\gamma$ and using the Cauchy-Schwarz inequality for sequences. This yields
\begin{equation}
\iint_{\R^3\times\R^3}\frac{|\gamma(x,y)|^2}{|x-y|}dx\,dy\leq D(\rho_\gamma,\rho_\gamma)<\ii,
\end{equation}
showing that the exchange term is also a well-defined  quantity when $(\gamma,\alpha)\in\cK$.

To deal with the pairing term in $\cE(\gamma,\alpha)$, we recall the Hardy-Kato inequality 
\begin{equation}
\frac1{|x|}\leq \frac\pi2 \sqrt{-\Delta} \, ,
\label{Kato}
\end{equation}
see \cite{Kato,Herbst-77}. Applying \eqref{Kato} in the variable $x$ with $y$ fixed, we obtain
\begin{equation}
\iint_{\R^3\times\R^3}\frac{|\alpha(x,y)|^2}{|x-y|}dx\,dy\leq \frac\pi2\tr\big(\sqrt{-\Delta}\,\alpha\alpha^*\big)<\ii,
\label{estim_pairing_term1}
\end{equation}
for any $(\gamma, \alpha) \in \cK$. Next, we introduce the $2\times2$-matrix
$$
\Gamma=\left(\begin{matrix}
\gamma&\alpha\\ \alpha^*&1-\overline{\gamma}\end{matrix}\right)$$ 
which defines an operator acting on $L^2(\R^3; \C^q)\oplus L^2(\R^3; \C^q)$.
Following the convention in \cite{BacLieSol-94}, we refer to any such $\Gamma$ satisfying the constraint \eqref{constraint_2} as an {\em admissible 1-particle density matrix (1-pdm)}. In particular, we immediately obtain from $0\leq\Gamma\leq1$ that $\Gamma^2\leq\Gamma$, and this leads to the operator inequality 
\begin{equation}
\gamma^2+\alpha\alpha^*\leq \gamma ,
\label{estim_pairing}
\end{equation}
as already observed in \cite{BacLieSol-94}. This fact combined with \eqref{estim_pairing_term1} yields 
\begin{equation}
\iint_{\R^3\times\R^3}\frac{|\alpha(x,y)|^2}{|x-y|} dx\,dy\leq \frac\pi2\tr\big(\sqrt{-\Delta}\,\gamma\big) .
\label{estim_pairing_term}
\end{equation}
Hence the pairing term is controlled by the kinetic energy.

\begin{remark}
Because of \eqref{estim_pairing_term}, we will have to assume that $0\leq\kappa< 4/\pi$ in the following; see also Remark \ref{rem:kappa2} below why this upper bound on $\kappa$ is optimal. With regard to its physical application to stars, such a smallness condition for $\kappa$ is harmless, since we have $\kappa = G m^2 \sim 10^{-38}$, when $G$ is Newton's gravitational constant and $m$ the neutron mass. Moreover, it is also natural to consider the HFB model in the limiting regime when $\kappa\ll1$ with $\kappa\lambda^{2/3}$ being fixed (i.\,e., the semi-classical regime for relativistic fermions).  \label{rem:kappa}
\end{remark}

\begin{remark}
 The arguments presented above can be easily adapted to prove that $\cE : \cK \to \R$ is continuous with respect to the norm $\| \cdot \|_\cX$.
\end{remark}

All this proves that the HFB energy functional $(\gamma,\alpha)\mapsto\cE(\gamma,\alpha)$ with mass parameter $m$ and coupling constant $\kappa$
is a well-defined functional on $\cK$. Hence it makes sense to consider the following minimization problem 
\begin{equation} \label{eq:IHF}
\boxed{ I(\lambda) = \inf  \big \{ \cE(\gamma,\alpha) : (\gamma,\alpha) \in \cK_{\lambda} \big \} }
\end{equation}
with $\lambda \geq 0$ given. Note we have $I(0)=0$ with unique minimizer $(\gamma,\alpha)=(0,0)$.

Before stating our main results on the variational problem (\ref{eq:IHF}), we need to define the largest possible number $\lambda$ of particles of the system for which the energy is bounded below. We have the following statement. 

\begin{prop}{\bf (Chandrasekhar Limit).} \label{prop:def_max_mass}
Let $m \geq 0$ and $0\leq\kappa<4/\pi$ be given. Then there exists a unique number $\lambda^{\rm HFB}(\kappa)>0$, which is independent of $m$, such that the following holds.

\begin{enumerate}
\item[$(i)$] For $0\leq\lambda\leq \lambda^{\rm HFB}(\kappa)$, we have $I(\lambda)>-\ii$.
\item[$(ii)$] For $\lambda>\lambda^{\rm HFB}(\kappa)$, we have $I(\lambda)=-\ii$.
\end{enumerate}

Furthermore, the function $\lambda^{\rm HFB}(\kappa)$ is nonincreasing and continuous with respect to $\kappa$. It satisfies the asymptotic estimate
\begin{equation}
\lambda^{\rm HFB}(\kappa)\sim\left ( \frac{\tau_c}{\kappa} \right )^{3/2}q^{-1/2} \quad \mbox{as} \quad \kappa \to 0,
\label{behavior_critical_mass}
\end{equation}
for some universal constant $\tau_c\simeq2.677$.
\end{prop}

The proof of Proposition \ref{prop:def_max_mass} is given in Appendix \ref{appendix:proof_def_max_mass} below. The constant $\tau_c$ was already defined in \cite{LieYau-87}, where it appears from a Thomas-Fermi-type energy obtained in the semi-classical limit $\kappa\to0$ with $\kappa\lambda^{2/3}=O(1)$.
 
\begin{remark} \label{rem:kappa2} 
It will be a consequence of our proof that  $I(\lambda)=-\ii$ for all $\lambda>0$ and all $m\geq0$ when $\kappa > 4/\pi$; see Remark \ref{rmk:unbounded} below. Saying differently, we have $\lambda^{\rm HFB}(\kappa)=0$ when $\kappa > 4/\pi$. 
\end{remark}

\begin{remark}
One simplified model consists in restricting the variational set to HFB states without pairing, i.\,e., states of the form $(\gamma,0)$ with $0\leq\gamma\leq 1$ and $\tr(\gamma)=\lambda$. A further simplification consists in restricting the energy to projectors only, i.\,e., $\gamma^2=\gamma$, when $\lambda=N$ is an integer (this is the usual Hartree-Fock model). The so-obtained functionals satisfy similar properties as $\cE(\gamma, \alpha)$ on the corresponding variational sets, except that the assumption $\kappa<4/\pi$ can be removed. 
\end{remark}

\begin{remark}
It can be easily seen from our proof provided in Appendix \ref{appendix:proof_def_max_mass}, that for every $0\leq \tau<\tau_c$
$$\lim_{\substack{\lambda\to\ii\\ \kappa \lambda^{2/3}\to \tau}}\frac{I(\lambda)}{I_{\rm Ch}(\lambda,\kappa)}=1$$
where $I_{\rm Ch}(\lambda,\kappa)$ is the Chandrasekhar ground state energy as defined in \cite{LieThi-84,LieYau-87}.
\end{remark}

\section{Main Results} \label{sec:main}

In this section, we state our main results for the HFB variational problem defined in (\ref{eq:IHF}). First, we establish the existence of minimizers in Theorem \ref{thm:relaxed}. Essential properties of minimizers will then be formulated in Theorem \ref{thm:prop_min} for the general case, and in Theorem \ref{thm:infinite_rank} for the physically relevant $q=2$ case. Furthermore, we study the simplified Hartree-Fock model and derive an existence result of minimizers; see Theorem \ref{thm:orbitals} below. Finally, we briefly discuss the time-dependent HFB theory.

\subsection{Existence of Minimizers in HFB Theory}\label{sec:exists}
We first formulate the existence result of minimizers with particle numbers below the critical threshold. Since the functional $\cE$ is translation invariant, we have to take into account the unitary action $\tau_y$ on $L^2(\R^3)$ given by the group of translations in $\R^3$; see Equation (\ref{eq:tau}) above for the definition of $\tau$. The precise existence result now reads as follows.

\begin{thm}{\bf (Existence of Minimizers).} \label{thm:relaxed} Fix the integer $q \geq 1$ describing the internal spin degrees of freedom. Furthermore, suppose that $m>0$ and $0<\kappa<4/\pi$. Then, for all $0 < \lambda <  \lambda^{\rm HFB}(\kappa)$, the following properties hold.

\begin{enumerate}

\item[$(i)$] Every minimizing sequence $\{ (\gamma_n,\alpha_n) \}_{n \in \N}$ for $I(\lambda)$ is relatively compact in $\cX$ up to translations. That is, there is a sequence $\{ y_n \}_{n \in \N} \subset \R^3$ such that, after passing to a suitable subsequence, we have
\[
\mbox{$\tau_{y_n}^* (\gamma_n,\alpha_n) \tau_{y_n} \rightarrow (\gamma,\alpha)$ strongly in $\cX$ as $n \to \infty$,}
\]
where $(\gamma,\alpha) \in \cK_\lambda$ is a minimizer for $I(\lambda)$. In particular, there exists a minimizer $(\gamma, \alpha)$ for $I(\lambda)$.

\item[$(ii)$] The following  binding inequality holds for all $0 < \lambda' < \lambda$:
\[
I(\lambda) < I(\lambda - \lambda') + I(\lambda'), 
\] 
\end{enumerate}
\end{thm}

\noindent
Let us make some remarks.

\begin{remark}\label{rmk:general_potential}
As can be checked from our proof, Theorem \ref{thm:relaxed} also holds true when the purely gravitational interaction $-1/|x-y|$ is replaced by a radial potential $W(|x-y|)$ satisfying the following assumptions
\begin{equation*}
\label{ass:W} 
|W(|x|)|\leq \frac{1}{|x|}\quad\text{for all $x\in\R^3$,}\qquad\text{and}\qquad W(|x|)\leq -\frac{\epsilon}{|x|}\quad\text{for $|x|\geq R_0$,}
\end{equation*}
for some constants $\epsilon >0$ and $R_0 > 0$. Of course, the largest particle number $\lambda^{\rm HFB}(\kappa)$ has to be defined accordingly.
\end{remark}

\begin{remark}
In fact, Properties (i) and (ii) are found to be equivalent, as it is well-known in the general framework of concentration-compactness method \cite{Lions-84} for variational problems with translation invariance.
\end{remark}

\begin{remark}
Note that we could not prove the existence of minimizers for $\lambda=\lambda^{\rm HFB}(\kappa)$, for which the energy is nevertheless bounded from below. The reason is that there seems to be a lack of coercivity in this case (see Lemma \ref{lem:coercive} below). Furthermore, if $m > 0$ holds, we do no expect minimizers to exist with critical mass $\lambda = \lambda^{\rm HFB}(\kappa)$. See also \cite{LieYau-87} for a related nonexistence result for a scalar Hartree equation describing boson stars.

\end{remark}

\begin{remark}
It is essential to assume that $m >0$ holds, which is also the physically relevant scenario. In the case of vanishing mass $m=0$, a further loss of compactness (besides the one induced by translation invariance) arises from scaling. 
\end{remark}

\subsection{Properties of HFB Minimizers}

\subsubsection{Nonlinear Equation and Decay Estimate}
Let us now describe some properties of HFB minimizers $(\gamma,\alpha)$ obtained in Theorem \ref{thm:relaxed}. As before, we denote
\begin{equation}
\Gamma=\left(\begin{matrix}
\gamma&\alpha\\ \alpha^*&1-\overline{\gamma}\end{matrix}\right), 
\label{def_GAMMA}
\end{equation}
and we introduce the following \emph{HFB mean-field} operator 
\begin{equation}
F_{\Gamma}:= \left(\begin{matrix}
H_\gamma & -\kappa\frac{\alpha(x,y)}{|x-y|}\\
-\kappa\frac{\alpha^*(x,y)}{|x-y|} & -\overline{H_\gamma}\\
\end{matrix}\right)
\label{def_H_SCF2}
\end{equation}
acting on $L^2(\R^3; \C^q)\oplus L^2(\R^3; \C^q)$. Here
\begin{equation}
H_\gamma:= T-\kappa(\rho_\gamma\ast\frac{1}{|\cdot|})(x)+\kappa\frac{\gamma(x,y)}{|x-y|}
\label{def_H_SCF}
\end{equation}
is the usual mean-field operator of Hartree-Fock theory \cite{LieSim-77,Lions-87,BacLieSol-94}. Moreover, it turns to be convenient to define 
\begin{equation}
N=\left(\begin{matrix}
1&0\\ 0&-1\end{matrix}\right).
\label{def_N}
\end{equation}
Note that $\Gamma$ commutes with $N$ if and only if $\alpha=0$, that is if and only if the corresponding quasi-free state in the Fock space also commutes with the number operator $\cN$; see, e.\,g., \cite{BacLieSol-94} for the Fock space formalism in HFB theory.

We can now state some fundamental properties of the minimizers for the HFB model considered in this paper.  

\begin{thm}{\bf (Properties of Minimizers).} \label{thm:prop_min} Let $q \geq 1$ be given and suppose $m>0$ and $0<\kappa < 4/\pi$. Assume that $(\gamma,\alpha)\in\cK_\lambda$ is a minimizer for $I(\lambda)$ for some $\lambda >0$. Then there exists a negative real number $\mu<0$ such that $\Gamma=\Gamma(\gamma, \alpha)$ solves the following nonlinear equation
 \begin{equation}
\Gamma=\chi_{(-\ii,0)}\left(F_\Gamma-\mu N\right)+D ,
\label{SCF_GAMMA} 
\end{equation}
where $D$ is a finite rank operator of the same matrix form as $\Gamma$ and satisfies ${\rm ran} (D)\subset \ker (F_\Gamma-\mu N)$. If moreover $\alpha \neq 0$ holds, we have the improved bound $\mu < \beta - m$, where the number  $\beta < m$ is given by Proposition \ref{prop:comparison_I_J} below.

Finally, the following decay estimate for the density function holds, for all $R >0$ sufficiently large,
$$
\int_{|x| \geq R} \rho_\gamma(x) \, dx \leq \frac{C}{R^2}, 
$$
where $C > 0$ is some constant independent of $R$. 
\end{thm}

\begin{remark}
Again the same result holds true when the purely attractive Newtonian interaction is replaced by a potential $W$ satisfying the assumptions in Remark \ref{ass:W}.
\end{remark}

\begin{remark}
It is an interesting open question whether the HFB-minimizers $(\gamma, \alpha) \in \cK_\lambda$ (or less ambitiously) the density $\rho_\gamma$ are unique for given $\tr(\gamma) = \lambda$, modulo symmetries of the HFB energy functional.
\end{remark}

The decay estimate for $\rho_\gamma(x)$ and a slightly refined version (see Lemma \ref{lem:rho_decay} below) will play an important role in the proof of Theorem \ref{thm:relaxed}. Let us note that if $\gamma$ were a finite-rank operator, the self-consistent equation (written in the same form as in \cite{Lewin-04a}) would show immediately that the associated density $\rho_\gamma(x)$ is indeed exponentially decaying. However, we believe that $\gamma$ is in general \emph{not a finite-rank operator} for any $q \geq 1$. We shall justify this claim for the physically relevant case $q=2$ below, provided we assume that $\alpha \neq 0$ holds.

\subsubsection{The spin-$1/2$ case ($q=2$)}
Let us now consider the specific case when the number of internal spin degrees of freedom is $q=2$, which corresponds to the physically relevant case of spin-$1/2$ fermions such as neutrons. In this case, it is shown by Bach-Fr\"ohlich-Jonsson in \cite{BacFroJon-08}, based on a concavity result in Lieb \cite{Lieb-73}, that the following equality holds:
$$
\boxed{ \text{$I(\lambda) = J(\lambda)$ for $q=2$ spin degrees of freedom.} }
$$
Here
\begin{equation}
 J(\lambda)=\inf\big\{2\cF(\tau)\ :\ \tau\in\cB\big(L^2(\R^3;\C)\big),\ \tau=\tau^*=\overline{\tau},\  \tr(\tau)=\lambda/2\big\},
\label{def_J}
\end{equation}
and the reduced (no-spin) energy is defined as
\begin{multline}
\cF(\tau):=\tr \left ( T \tau \right ) - \kappa D(\rho_{\tau}, \rho_{\tau} ) \\ + \frac{\kappa}{2} \Ex(\tau) 
-\frac{\kappa}{2} \iint_{\R^3\times\R^3}\frac{\left|\sqrt{\tau(1-\tau)}(x,y)\right|^2}{|x-y|}dx\,dy.
\label{def_reduced_energy}
\end{multline}
Furthermore, all the minimizers for $\cE(\gamma, \alpha)$ and $F(\tau)$, respectively, satisfy 
\begin{equation}
\gamma=\tau\otimes\left(\begin{matrix}1&0\\
0&1\end{matrix}\right)\quad\text{ and }\quad \alpha=\pm\sqrt{\tau(1-\tau)}\otimes\left(\begin{matrix}0&1\\
-1&0\end{matrix}\right).
\label{eq:tau_min}
\end{equation}

We note that $\cF(\tau)$ has a form similar to the M\"uller functional which was studied in \cite{FraLieSeiSie-07}. Indeed, a straightforward adaptation of an argument given in \cite{FraLieSeiSie-07} combined with (\ref{eq:tau_min}) leads to the following statement (its proof is outlined in Section \ref{sec:proof_infinite_rank}).

\begin{thm}{\bf (Infinite Rank of Minimizers if $\alpha \neq 0$).} Assume that $(\gamma, \alpha)$ is a minimizers for $I(\lambda)$ with $0< \kappa < 4/\pi$ and $\lambda > 0$. Then, if $\alpha \neq 0$ holds, the operators $\gamma$ and $\alpha$ both have infinite rank. 
\label{thm:infinite_rank}
\end{thm}

\begin{remark}
Since it is expected that $\alpha \neq 0$ holds in general, this result justifies the variational arguments developed below to prove a decay estimate for $\rho_\gamma$.
\end{remark}

\subsection{Minimizers for Simplified HF Model} \label{sec:HF_theory}
As we have mentioned above, a rough approximation consists in neglecting the pairing term, i.\,e., by restricting to density matrices having $\alpha=0$. A result similar to Theorem \ref{thm:relaxed} holds true for an adequately defined critical mass $\lambda^{\alpha=0}(\kappa)$ and without the assumption that $\kappa<4/\pi$. Minimizers then satisfy the following equation $\gamma=\chi_{(-\ii,\mu)}(H_\gamma)+\delta$. The proof is indeed much simpler, as can be seen from the proof of Theorem \ref{thm:relaxed}.

In the Hartree-Fock (HF) model, one restricts the search of a minimum to the variational set to density matrices $\gamma$ that are projectors of rank $N$, i.\,e. $\gamma=\sum_{i=1}^N|\phi_i\rangle\langle\phi_i|$,
with some $L^2$-orthonormal functions $\{ \phi_i \}_{i=1}^N$. For any integer $N \geq 1$, it is convenient to define the set
\begin{equation}
\cP_N:= \{(\gamma,0)\in\cK\ :\ \gamma^2=\gamma \},
\label{eq:cP_N}
\end{equation}
and accordingly the HF minimization is given by
\begin{equation}
\boxed{ \IHF(N):=\inf \big \{\cE(\gamma,0) : (\gamma,0)\in\cP_N \big \} .}
\label{def_HF_min}
\end{equation}
In analogy to the HFB model discussed above, we denote by $N^{\rm HF}(\kappa)$ the largest integer for which $\IHF(N)>-\ii$ holds. 
We now have the following result about the HF minimization problem.

\begin{thm} {\bf (Existence of Minimizers for HF-Model).} \label{thm:orbitals} Let $m>0$ and $\kappa>0$. 
First, there is no minimizer for $\IHF(1)$.
Then for every integer $N$ such that
$1 < N < N^{\rm HF}(\kappa)$,
the following holds:

\smallskip

\noindent$(i)$ Every minimizing sequence $\{ (\gamma_{n},0) \}_{n \in \N}$ for $\IHF(N)$ is relatively compact in $\cX$ up to translations. In particular, there is a minimizer $(\gamma,0) \in\cP_N$ for $\IHF(N)$.

\smallskip

\noindent$(ii)$ The following binding inequality holds, for all integers $1 \leq k \leq N-1$:
\[
\IHF(N) < \IHF(N-k) + \IHF(k).
\] 

\smallskip

\noindent$(iii)$ Any minimizer $(\gamma,0) \in \cP_N$ for $\IHF(N)$ can be written as $\gamma = \sum_{n=1}^N |\varphi_N \rangle \langle \varphi_N |$, where the functions $\varphi_n \in \Hhalf$ satisfy 
\begin{equation}
H_\gamma \varphi_n = \lambda_n \varphi_n ,
\label{eq_SCF_phi_thm} 
\end{equation}
with $\lambda_1 \leq \cdots \leq \lambda_N < 0$ being $N$ negative eigenvalues of the the mean-field operator $H_\gamma$ defined in \eqref{def_H_SCF}.
\end{thm}

Since the HF model introduced above is an $N$-body theory, one can prove, by using the geometrical methods of \cite{Enss-77,Simon-77,Sigal-82}, that every minimizing sequence is compact (up to translation) if and only if the quantized binding inequality in $(ii)$ holds true. This was done for the first time by Friesecke in \cite{Friesecke-03} in the context of the usual atomic Hartree-Fock and MCSCF theories. The proof of Theorem \ref{thm:orbitals} relies on Friesecke's method, together with some ideas of the proof of Theorem \ref{thm:relaxed}. In some sense it is much easier than for Theorem \ref{thm:relaxed}, since it is possible to prove the quantized binding inequalities by induction on $N$. We briefly explain this in Section \ref{sec:proof_thm_orbitals} below.

\subsection{Time-Dependent Theory}\label{sec:time_dependent}

We now make a brief digression into the related {\em time-dependent HFB theory.} As it is well-known from nonlinear Schr\"odinger equations, a detailed understanding of the variational calculus for the time-independent theory turns out to be of great use when addressing dynamical questions; e.\,g., the stability of solitary waves and blowup analysis. Here, however, we shall mostly limit ourselves to the dynamical (orbital) stability of solutions generated by minimizers (which essentially follows from Theorem \ref{thm:relaxed} as a by-product). As a concluding outlook, we briefly comment on the finite-time blowup for the time-dependent HFB equations, corresponding to the onset of the dynamical collapse of a neutron star.      

For the reader's orientation, we recall that the time-dependent HFB equations generated by the functional $\cE(\gamma,\alpha)$ can be written in commutator form as 
\begin{equation} \label{eq:Neumann}
i \partial_t \Gamma = [F_\Gamma, \Gamma ] .
\end{equation}  
Here and as usual,  $\Gamma = \Gamma(\gamma, \alpha)$ denotes the admissible 1-pdm defined in \eqref{def_GAMMA} and $F_\Gamma$ is the HFB mean-field operator introduced in \eqref{def_H_SCF2}. Clearly, the equation \eqref{eq:Neumann} is a nonlinear evolution equation. With regard to the wellposedness of its initial-value problem, we note that a straightforward adaptation of \cite{ChaGla-75,Chadam-76,BovPraFan-76,FroLen-07,HaiLewSpa-05,Lenzmann-07} yields the following results.

\medskip

\paragraph{\bf Local Wellposedness.} For each initial datum $(\gamma_0, \alpha_0) \in \cK$ and $T>0$ sufficiently small, there exists a unique solution $(\gamma,\alpha) \in C^0([0,T); \cX)\cap C^1([0,T); \cX')$ solving \eqref{eq:Neumann} such that $\Gamma(0) = \Gamma(\gamma_0,\alpha_0)$. Moreover, we have conservation of  energy and expected number of particles, i.\,e.,
\[
\forall 0 \leq t < T,\qquad \cE(\gamma(t),\alpha(t)) = \cE(\gamma_0, \alpha_0) \quad \mbox{and} \quad \tr (\gamma(t)) = \tr(\gamma_0).
\]

\medskip

\paragraph{\bf Global Wellposedness.} If the initial datum $(\gamma_0, \alpha_0) \in \cK$ satisfies the smallness condition $\tr(\gamma_0) < \lambda^{\rm HFB}(\kappa)$, then we have the a priori bound
\[
\sup_{0 \leq t < T} \| (\gamma(t), \alpha(t)) \|_{\cX} \leq C
\]
and the solution $\Gamma(t) = \Gamma(\gamma(t), \alpha(t))$ extends to all times $0 \leq t < \infty$.

\medskip

Next, we consider the behavior of minimizers for the HFB energy $\cE(\gamma,\alpha)$ with respect to the evolution equation \eqref{eq:Neumann}. To this end, let us assume that the initial condition $(\gamma_0, \alpha_0)$ is a minimizer for $I(\lambda)$ as given by Theorem \ref{thm:prop_min}. With this choice of initial conditions and using equation \eqref{SCF_GAMMA} above, an elementary calculation shows that the corresponding solution $\Gamma(t) = \Gamma(\gamma(t), \alpha(t))$ of \eqref{eq:Neumann} is given by
\begin{equation}
\gamma(t)  =\gamma_0 \quad \mbox{and} \quad \alpha(t) = e^{-2\mu it} \alpha_0 ,
\end{equation}
or equivalently in a $2\times2$ matrix form
\begin{equation}
\Gamma(t)=e^{-i\mu Nt}\,\Gamma_0\,e^{i\mu Nt}
\end{equation}
with $\mu <0$ taken from Theorem \ref{thm:prop_min} and $N$ defined in \eqref{def_N}. Hence HFB minimizers give rise to stationary solutions of \eqref{eq:Neumann}, as one naturally expects. Moreover, by adapting a well-known general argument in \cite{CazLio-82}, the relative compactness result of Theorem \ref{thm:relaxed} and the conservation laws imply orbital stability of minimizers. The precise statement can be formulated as follows.
Let us define the set of minimizers
\[
\mathcal{M}_\lambda = \big \{ (\gamma, \alpha) \in \cK_{\lambda} : \cE(\gamma,\alpha) = I(\lambda ) \big \},
\] 
and introduce the distance function  on $\cK$
$${\rm dist}_{\mathcal{M}_\lambda}(\gamma_0, \alpha_0) = \inf_{(\gamma, \alpha) \in \mathcal{M}_\lambda} \| (\gamma_0, \alpha_0) - (\gamma, \alpha) \|_{\cX}.$$

\begin{thm} {\bf (Stability of Minimizers under HFB Time Evolution).} Under the assumptions of Theorem \ref{thm:relaxed}, the set of minimizers $\cM_\lambda$ for $I(\lambda)$ with $0 < \lambda < \lambda^{\rm HFB}(\kappa)$ is {\em orbitally stable} in the following sense: for every $\epsilon >0$, there is $\delta > 0$ such that if the initial condition satisfies
${\rm dist}_{\mathcal{M}_\lambda}(\gamma_0, \alpha_0) < \delta$, then the corresponding solution $\Gamma(t) = (\gamma(t),\alpha(t))$ of (\ref{eq:Neumann}) exists for all times $t \geq 0$ and obeys $\sup_{t \geq 0} {\rm dist}_{\mathcal{M}_\lambda}(\gamma(t), \alpha(t)) < \epsilon$.
\end{thm}

Finally, we briefly comment on finite-time blowup for the time-dependent HFB equation. Indeed, recent work by Hainzl and Schlein \cite{HaiSch-09} (which extends previous results in \cite{FroLen-07, FroLen-07B}) addresses blowup for the Hartree-Fock (HF) approximation. Since, in particular, the HFB time evolution (\ref{eq:Neumann}) preserves HF states, we can immediately infer from \cite{HaiSch-09, FroLen-07} the following result.

\medskip

\paragraph{\bf Finite-Time Blowup for HF States.} Let the initial datum $(\gamma_0, 0) \in \mathcal{P}_N$ be a Hartree-Fock state; see \eqref{eq:cP_N}. Furthermore, assume that $\gamma_0 \in C^\infty_0(\R^3 \times \R^3)$ and that we have {\em radially symmetry} in the sense that $\gamma_0(x,y) = \gamma_0(Rx,Ry)$ for all $R \in SO(3)$. Then, the energy condition
$\cE(\gamma_0,0) < -m N$
implies that the corresponding solution $\Gamma(t) = \Gamma(\gamma(t), 0)$ of \eqref{eq:Neumann} blows up in finite time, i.\,e., we have
$\lim_{t \nearrow T} \tr (\sqrt{-\Delta} \gamma(t)) = +\infty,$
for some $0 < T < \infty$. Furthermore, the solution exhibits minimal mass concentration at the origin in the sense that, for every $R > 0$, 
\[
\liminf_{t \nearrow T} \int_{|x| < R} \rho_{\gamma(t)} \, dx \geq N^{\rm HF}(\kappa) ,
\]
where $N^{\rm HF}(\kappa) > 0$ is the same constant as in Theorem \ref{thm:orbitals} above.

\begin{remark}
In order for $\cE(\gamma_0, 0) < -mN$ to hold, we must have large initial data in the sense that $\tr(\gamma_0) = N > N^{\rm HF}(\kappa)$ holds (which implies that $I^{\rm HF}(N) =-\infty$ for such large $N$); cf.~Section \ref{sec:HF_theory} for the definition of $N^{\rm HF}$. 
\end{remark}

\begin{remark}
It is of interest to extend the above blowup result to HFB states with non-vanishing pairing $\alpha \neq 0$; or, even more ambitious, to relax the radiality assumptions on the initial data.
\end{remark}


\section{Preliminaries} \label{sec:prelim}

In this section we collect some preliminary results that are needed to set up the variational calculus for the proof of the main theorems.

\subsection{Weak-$*$ Topology on $\cX$} \label{subsec:weaktopo}
Recall that we assume that $q \geq 1$ (describing the spin degrees of freedom) is a fixed integer and that $\gS_p$ with $1 \leq p < \infty$ denote the Schatten spaces of bounded operator acting on $L^2(\R^3;\C^q)$. It is well-known that $\gS_p$ are reflexive when $1<p<\ii$, the dual being $\gS_{p'}$ where $p^{-1}+{p'}^{-1}=1$. Furthermore, it is a classical fact that $\gS_1$ is the dual of the space of compact operators $\gS_\ii$; see \cite{Simon-79}. This induces naturally a {\em weak-$*$ topology} on $\cX$ for which the unit ball is compact, by the Banach-Alaoglu theorem. More precisely, we say that $(\gamma_n,\alpha_n) \wto (\gamma,\alpha)$ weakly-$\ast$ in $\cX$ if
\begin{equation} \label{def:weak}
\tr \big ( (1-\Delta)^{1/4} \gamma_n (1-\Delta)^{1/4} K \big ) \rightarrow \tr \big ( (1-\Delta)^{1/4} \gamma (1-\Delta)^{1/4} K \big ), 
\end{equation}
\begin{equation}
\tr \big ( (1-\Delta)^{1/4} \alpha_n K_2 \big ) \rightarrow \tr \big ( (1-\Delta)^{1/4} \alpha K_2\big ), 
\end{equation}
for all $K \in \gS_\ii$ and all $K_2\in\gS_2$. Note that the convex set of fermionic density matrices $\cK \subset \cX$ is both closed in the strong topology of $\cX$ and in the weak-$*$ topology. However, the set $\cK_\lambda$ is {\em not closed} in the weak-$*$ topology. That is, we can have $(\gamma_n,\alpha_n) \wto (\gamma,\alpha)$ weakly-$*$ in $\cX$ with $\gamma_n \in \cK_\lambda$ such that $\tr(\gamma) < \liminf \tr(\gamma_n)$. Put differently, the continuous linear functional $(\gamma,\alpha) \mapsto \tr(\gamma)$ on $\cX$ is not weakly-$*$ continuous. 

Throughout this paper, we shall make use of the following important fact that helps conclude strong convergence a posteriori: Suppose $(\gamma_n,\alpha_n) \wto (\gamma,\alpha)$ weakly-$*$ in $\cX$ and moreover assume that
\begin{equation}
\tr \big ( (1-\Delta)^{1/4} \gamma_n (1-\Delta)^{1/4} \big ) \rightarrow \tr \big ( (1-\Delta)^{1/4} \gamma (1-\Delta)^{1/4} \big ), 
\end{equation}
and
\begin{equation}
\tr \big ( (1-\Delta)^{1/4} \alpha^*_n\alpha_n (1-\Delta)^{1/4} \big ) \rightarrow \tr \big ( (1-\Delta)^{1/4} \alpha^*\alpha(1-\Delta)^{1/4} \big ) .
\end{equation}
Then we have $(\gamma_n,\alpha_n) \to (\gamma,\alpha)$ strongly in $\cX$. This follows from \cite[Thm A.6]{Simon-79}.

\subsection{Some Basic Properties of HFB Energy}

We have already seen that the energy $\cE$ is well-defined on $\cK$. Next, we collect some basic facts.

\begin{lemma}{\bf (Basic Properties of $I(\lambda)$).} \label{lem:IHF} Let $m>0$ and $0<\kappa<4/\pi$. Then the following holds.

\smallskip

\noindent$(i)$ For all $0 < \lambda < \lambda^{\rm HFB}(\kappa)$, we have $-\infty < I(\lambda) < 0$ and, for all $0 < \lambda' < \lambda$,
\begin{equation*}
 I(\lambda) \leq I(\lambda-\lambda') + I(\lambda'). 
\end{equation*}
Moreover, if equality holds for some $\lambda$ and $\lambda'$, then there exists a minimizing sequence for $I(\lambda)$ that is not relatively compact in $\cX$ up to translations.

\smallskip

\noindent$(ii)$ The function $\lambda \mapsto I(\lambda)$ is continuous and decreasing. 

\smallskip

\noindent$(iii)$ The functional $\cE$ is \underline{not} weakly-$\ast$ lower semicontinuous on $\cK_\lambda$.
\end{lemma}
\begin{proof}
The proof of the subadditivity inequality in Part $(i)$ follows from classical arguments: we take two states $(\gamma,\alpha)\in\cK_{\lambda-\lambda'}$ and $(\gamma',\alpha')\in\cK_{\lambda'}$ that are almost minimizers of $I(\lambda-\lambda')$ and $I(\lambda')$, respectively. By a density argument, we can assume that $(\gamma, \alpha)$ and $(\gamma',\alpha')$ are both finite-rank operators with compact support. Next, we define the sequence of states 
\begin{equation}
(\widetilde{\gamma}_n, \widetilde{\alpha}_n) =  (\gamma,\alpha)+\tau^*_{n e}(\gamma',\alpha')\tau_{n e}, \quad \mbox{with $n=1,2,3,\ldots$}
\end{equation}
where $e \in\R^3$ is some fixed unit vector and $\tau_{ne}f = f(\cdot - ne)$. Note that $(\gamma, \alpha)$ and $\tau^*_{ne} (\gamma',\tau') \tau_{ne}$ commute for $n$ large enough, and we easily verify that $(\widetilde{\gamma}_n, \widetilde{\alpha}_n) \in \cK_\lambda$ when $n$ is sufficiently large. By taking the limit $n \to \infty$, it is straightforward to conclude the subadditivity inequality. If moreover equality holds for some $\lambda'$, then $(\widetilde{\gamma}_n, \widetilde{\alpha}_n)$ furnishes a minimizing sequence that fails to be relatively compact in $\cX$ up to translations.

To see that $I(\lambda)<0$ holds, we choose a fixed state $(\gamma,0) \in \cK_\lambda$ with $\gamma$ smooth and finite-rank and such that $D(\rho_\gamma,\rho_\gamma)-\Ex(\gamma)>0$ holds. (The later condition is indeed equivalent to saying that ${\rm rank} \, (\gamma) \geq 2$.) Next, we consider the energy of the rescaled state $\gamma_\delta:=U_\delta\gamma U_\delta^*$ where $(U_\delta f)(x)=\delta^{-3/2}f(x/\delta)$. Note that $(\gamma_\delta, 0) \in \cK_\lambda$ for all $\delta > 0$, as one easily verifies. Using now the operator inequality $T \leq -\frac{1}{2m} \Delta$, we deduce that
\begin{equation*}
\cE(\gamma_\delta,0) \leq \frac{1}{\delta^2} \tr ( -\frac{1}{2m} \Delta \gamma )  - \frac{\kappa}{2 \delta}  (D(\rho_\gamma,\rho_\gamma)-\Ex(\gamma))<0, 
\end{equation*}
provided that $\delta > 0$ is chosen sufficiently large. This establishes Part $(i)$.

To show Part $(ii)$, we note that the strict monotonicity of $I(\cdot)$ is an obvious consequence of the subadditivity inequality and the strict negativity of $I(\cdot)$. The continuity of $I(\cdot)$ follows readily from using trial states. This proves Part $(ii)$.

To prove Part $(iii)$, we fix some $0<\lambda^1<\lambda$ and consider a state $(\gamma^1,0) \in\cK_{\lambda^1}$ with compact support (meaning that its density $\rho_{\gamma^1}$ has a compact support). Furthermore, let $(\gamma^2,0) \in\cK_{\lambda-\lambda^1}$ be with compact support and such that $\cE(\gamma^2,0)<0$, which is possible because of $I_\alpha(\lambda-\lambda^1)<0$. Define then $\gamma_n:=\gamma^1+\tau_{ne}\gamma^2\tau_{ne}^*$ with $\tau_{ne}$ as above. For $n$ large enough, the terms in $\gamma_n$ have disjoint supports and we obtain that $(\gamma_n,0)\in\cK_\lambda$. Next, a computation shows that 
$$\lim_{n\to\ii}\cE(\gamma_n,0)=\cE(\gamma^1,0)+\cE(\gamma^2,0)<\cE(\gamma^1,0).$$
 Since $\gamma_n \wto \gamma^1$ weakly-$*$ in $\cX$, this strict inequality shows that $\cE(\gamma,\alpha)$ is not weakly-$*$ lower semicontinuous on $\cX$. \end{proof}

\begin{lemma} {\bf (Coercivity of $\cE$).} \label{lem:coercive} Let $m>0$ and $0\leq\kappa<4/\pi$. For $0 \leq \lambda < \lambda^{\rm HFB}(\kappa)$, the energy $\cE$ is coercive on $\cK_\lambda$, i.\,e., we have 
$$\cE(\gamma,\alpha) \rightarrow \infty \quad \mbox{as} \quad \|(\gamma,\alpha) \|_\cX \rightarrow \infty. $$ 
In particular, all minimizing sequences for $\cE$ on $\cK_\lambda$ with $0 \leq \lambda < \lambda^{\rm HFB}(\kappa)$ are bounded in $\| \cdot \|_\cX$-norm. 
\end{lemma}

\begin{proof}
For any $(\gamma, \alpha) \in \cK_\lambda$ and $1 > \epsilon > 0$, we have
\begin{equation} \label{ineq:coerc}
\cE(\gamma, \alpha) \geq \epsilon \tr (T \gamma) + (1-\epsilon) I_{\kappa/(1-\epsilon)}(\lambda),
\end{equation}
where $I_{\kappa/(1-\epsilon)}(\lambda)$ is the ground state energy of $\cE$ on $\cK_\lambda$ with $\kappa$ replaced by $\kappa/(1-\epsilon)$. By continuity of $\lambda^{\rm HFB}(\kappa)$ and since $0 \leq \lambda < \lambda^{\rm HFB}(\kappa)$, we can choose $\epsilon > 0$ small such that $\lambda \leq \lambda^{\rm HFB}(\frac{\kappa}{1-\epsilon})$ and hence $I_{\kappa/(1-\epsilon)}(\lambda) > -\infty$. 
Using that $\alpha^* \alpha \leq \gamma$ (see \eqref{estim_pairing}) and the strict positivity of $m >0$, we see that $\| (\gamma, \alpha) \|_{\cX} \to \infty$ implies that $\tr (T \gamma) \to \infty$. This fact, by (\ref{ineq:coerc}), implies that $\cE(\gamma, \alpha) \to \infty$ as well.\end{proof}

\subsection{Study of an Auxillary Functional}\label{sec:auxillary}
The HFB energy contains two non-convex terms: the direct and pairing term. It turns out to be expedient to introduce the following energy functional
\begin{equation}
\cG(\gamma,\alpha):=\tr(K\gamma)-\frac\kappa2\iint_{\R^3\times\R^3}\frac{|\alpha(x,y)|^2}{|x-y|}dx\,dy, 
 \label{def_F}
\end{equation}
for $(\gamma, \alpha) \in \cK$, where we set
\begin{equation}
K:=\sqrt{-\Delta + m^2}.
\end{equation}
Recall estimate \eqref{estim_pairing_term} which tells us that the pairing term is controlled by the kinetic energy. Hence we have the lower bound
\begin{equation}
\cG(\gamma,\alpha)\geq 0,
\end{equation}
when $0\leq \kappa\leq 4/\pi$. We start by stating a result similar to \cite[Theorem 1]{HaiLewSer-05b}.
\begin{prop} \label{prop:pairing_kinetic} Let $m>0$ and $0\leq\kappa<4/\pi$. 
Then $\cG$ is weakly-$\ast$ lower semicontinuous on $\cK$. This means if $\{(\gamma_n,\alpha_n)\}\subset\cK$ satisfies $(\gamma_n,\alpha_n)\wto (\gamma,\alpha)$ weakly-$\ast$ in $\cX$, then
\begin{equation}
\liminf_{n\to\ii}\cG(\gamma_n,\alpha_n)\geq \cG(\gamma,\alpha).
\label{F_wlsc} 
\end{equation}
Additionally, equality holds if and only if the convergence is strong in $\cX$.
\end{prop}
\begin{proof}
Consider a sequence in $\cK$ such that $(\gamma_n,\alpha_n)\wto(\gamma,\alpha)$ weakly-$\ast$ in $\cX$. By Lemma \ref{lem:estim_density}, $\sqrt{\rho_{\gamma_n}}$ is bounded in $H^{1/2}(\R^3)$. Thanks to a Rellich-type theorem for $\Hhalf$ and extracting a subsequence if necessary, we can assume that $\rho_{\gamma_n}\wto\rho_\gamma$ weakly in $L^p(\R^3)$ for $1 < p\leq 3/2$ and strongly in $L^p_{\rm loc}(\R^3)$ for $1\leq p< 3/2$. As for the sequence $\alpha_n(x,y)$, we notice that
\begin{equation*} 
\tr \big ( (1-\Delta)^{1/4} \alpha_n \alpha_n^* (1-\Delta)^{1/4} \big ) = \frac{1}{2} \big \langle \big ( \sqrt{1-\Delta_x} +\sqrt{1-\Delta_y} \big ) \alpha_n, \alpha_n \big \rangle_{L^2(\R^3 \times \R^3 )},
\end{equation*}
Thus $\alpha_n(x,y)$ is a bounded sequence in $H^{1/2}(\R^3\times\R^3)$, and we may assume that $\alpha_n(x,y)$ converges strongly in $L^2_{\rm loc}(\R^3\times\R^3)$, by a Rellich-type theorem again.

Next, we recall the definition of the cutoff functions $\chi_R$ and $\zeta_R$ in \eqref{eq:cutoff}. We then find the estimate
\begin{align*}
&\iint_{\R^3\times\R^3}\frac{|\alpha_n(x,y)|^2}{|x-y|}dx\,dy\nonumber\\
&\qquad= \iint_{\R^3\times\R^3}\frac{\zeta_{R}(x)^2|\alpha_n(x,y)|^2}{|x-y|}dx\,dy+\iint_{\R^3\times\R^3}\frac{\chi_{R}(x)^2|\alpha_n(x,y)|^2}{|x-y|}dx\,dy\nonumber\\
&\qquad\leq\frac\pi2\tr(K\zeta_{R}\gamma_n\zeta_{R})+\iint_{\R^3\times\R^3}\frac{\chi_R(x)^2\chi_{3R}(y)^2|\alpha_n(x,y)|^2}{|x-y|}dx\,dy+\frac{\norm{\alpha_n}^2_{\gS_2}}R.
\end{align*}
In the last line we have used  the Hardy-Kato inequality \eqref{Kato} and the inequality $\zeta_{R}\alpha_n\alpha_n^*\zeta_{R}\leq\zeta_{R}\gamma_n\zeta_{R}$ thanks to \eqref{estim_pairing}. Because $(\gamma_n,\alpha_n)$ is bounded in $\cX$, the last term in the above estimate is $O(R^{-1})$.

Now we localize the kinetic energy. Lemma \ref{lem:IMS1} tells us that the nonlocal operator $K$ satisfies
\begin{align}
\chi_R K\chi_R+\zeta_RK\zeta_R  & \leq K + \frac{1}{\pi} \int_0^\infty \frac{1}{s+K^2} \left (|\nabla \chi_R |^2+|\nabla \zeta_R |^2 \right ) \frac{1}{s+K^2} \, \sqrt{s} \, ds \nonumber \\
& \leq K + C/R^2,
\end{align}
using $\norm{\nabla\chi_R}^2_{L^\ii}+\norm{\nabla\zeta_R}^2_{L^\ii}\leq C/R^2$. Since $\tr(\gamma_n)$ is uniformly bounded, this gives us $\tr(K\gamma_n)\geq \tr(K\chi_{R}\gamma_n\chi_{R})+\tr(K\zeta_{R}\gamma_n\zeta_{R})-C/R^2$. Thus we obtain
\begin{multline}
\cG(\gamma_n,\alpha_n)\geq \tr\big(K(\chi_{R}\gamma_n\chi_{R}-\chi_{R}\alpha_n\chi_{3R}^2\alpha_n^*\chi_{R})\big)+\tr\big(K\chi_R\alpha_n\chi_{3R}^2\alpha_n^*\chi_R\big)\\
-\frac\kappa2\iint_{\R^3\times\R^3}\frac{\chi_R(x)^2\chi_{3R}(y)^2|\alpha_n(x,y)|^2}{|x-y|}dx\,dy +(1-\kappa\pi/4)\tr(K\zeta_R\gamma_n\zeta_R)-\frac{C}R.\label{estim_below_F_wlsc}
\end{multline}
Since $\alpha_n$ converges strongly to $\alpha$ in $L^2_{\rm loc}(\R^3\times\R^3)$, we have that $\chi_{R}\alpha_n\chi_{3R}^2\alpha_n^*\chi_{R}\to \chi_{R}\alpha\chi_{3R}^2\alpha^*\chi_{R}$ strongly in $\gS_1$. Hence we deduce that
$$\chi_{R}\gamma_n\chi_{R}-\chi_{R}\alpha_n\chi_{3R}^2\alpha_n^*\chi_{R}\wto \chi_{R}\gamma\chi_{R}-\chi_{R}\alpha\chi_{3R}^2\alpha^*\chi_{R}$$
weakly in $\gS_1$. Notice also that,  by \eqref{estim_pairing},
$$\chi_{R}\alpha_n\chi_{3R}^2\alpha_n^*\chi_{R}\leq \chi_{R}\alpha_n\alpha_n^*\chi_{R}\leq \chi_{R}\gamma_n\chi_{R}$$ 
Therefore we may use Fatou's Lemma to infer that 
$$\liminf_{n\to\ii}\tr\big(K(\chi_{R}\gamma_n\chi_{R}-\chi_{R}\alpha_n\chi_{3R}^2\alpha_n^*\chi_{R})\big)\geq \tr\big(K(\chi_{R}\gamma\chi_{R}-\chi_{R}\alpha\chi_{3R}^2\alpha^*\chi_{R})\big).$$
Next, we introduce the following sequence of functions on $\R^3 \times \R^3$ given by 
$\alpha_n^R(x,y):=\chi_R(x)\chi_{3R}(y)\alpha_n(x,y),$
which converges weakly in $H^{1/2}(\R^3\times\R^3)$ to $\alpha^R(x,y):=\chi_R(x)\chi_{3R}(y)\alpha(x,y)$. We can write
\begin{multline*}
 \tr\big(K\chi_R\alpha_n\chi_{3R}^2\alpha_n^*\chi_R\big)-\frac\kappa2\iint_{\R^3\times\R^3}\frac{\chi_R(x)^2\chi_{3R}(y)^2|\alpha_n(x,y)|^2}{|x-y|}dx\,dy\\
=\pscal{\left(K_x-\frac\kappa{2|x-y|}\right)\alpha^R_n,\alpha^R_n}_{L^2(\R^3\times\R^3)}.
\end{multline*}
Using that $\kappa<4/\pi$, the Hardy-Kato inequality yields
$$K_x-\frac\kappa{2|x-y|}\geq (1-\kappa\pi/4)K_x\geq (1-\kappa\pi/4)m>0 .$$
Thus we may use Fatou's Lemma one more time to obtain
\begin{multline}
\liminf_{n\to\ii}\pscal{\left(K_x-\frac\kappa{2|x-y|}\right)\alpha^R_n,\alpha^R_n}_{L^2(\R^3\times\R^3)}\\ \geq  \pscal{\left(K_x-\frac\kappa{2|x-y|}\right)\alpha^R,\alpha^R}_{L^2(\R^3\times\R^3)}.
\label{estim_alpha} 
\end{multline}
In summary, we have derived the following inequality
\begin{equation*}
\liminf_{n\to\ii}\cG(\gamma_n,\alpha_n)\geq \tr (K\chi_{R}\gamma\chi_{R}) -\frac\kappa2\iint_{\R^3\times\R^3}\frac{\chi_R(x)^2\chi_{3R}(y)^2|\alpha(x,y)|^2}{|x-y|}dx\,dy -\frac{C}R.
\end{equation*}
By taking the limit $R\to\ii$, we obtain the desired result 
$$\liminf_{n\to\ii}\cG(\gamma_n,\alpha_n)\geq \cG(\gamma,\alpha).$$

It remains to show the claim about strong convergence and equality. First, we note that if equality holds in \eqref{F_wlsc}, then by \eqref{estim_below_F_wlsc}, \eqref{estim_pairing} and $\kappa<4/\pi$,
\begin{equation}
 \lim_{R\to\ii}\limsup_{n\to\ii}\tr(K\zeta_R\gamma_n\zeta_R)=\lim_{R\to\ii}\limsup_{n\to\ii}\tr(K\zeta_R\alpha_n\alpha_n^*\zeta_R)=0 .
\label{lim_sup_R}
\end{equation}
On the other hand, the lower bound \eqref{estim_alpha} implies that in case of equality, for any fixed $R$, we have $\alpha^R_n\to\alpha^R$ strongly in $H^{1/2}(\R^3\times\R^3)$. We easily deduce from \eqref{lim_sup_R} that $\alpha_n\to\alpha$ strongly in $H^{1/2}(\R^3\times\R^3)$. By this strong convergence, we can now pass to the limit in the pairing term and we conclude
$\lim_{n\to\ii}\tr(K\gamma_n)=\tr(K\gamma)$, 
which itself eventually implies that $(\gamma_n,\alpha_n)\to(\gamma,\alpha)$ strongly in $\cX$. The proof of Proposition \ref{prop:pairing_kinetic} is now complete. \end{proof}

A simple but important consequence of Proposition \ref{prop:pairing_kinetic} is the following fact.

\begin{corollary}{\bf (Conservation of Mass implies Compactness).} \label{lem:compact} Let $m>0$ and $0<\kappa<4/\pi$. Consider a minimizing sequence $\{ (\gamma_n,\alpha_n) \}_{n \in \N}$ for $I(\lambda)$ with $\lambda < \lambda^{\rm HFB}(\kappa)$ and such that $(\gamma_n,\alpha_n) \wto (\gamma,\alpha)$ weakly-$*$ in $\cX$. Then $(\gamma_n,\alpha_n) \to (\gamma,\alpha)$ strongly in $\cX$ if and only if $\tr(\gamma) = \lambda$.
\end{corollary}

\begin{proof}
Clearly, if $\gamma_n \rightarrow \gamma$ strongly in $\cX$, then $\tr(\gamma) = \lambda$. To show that $\tr (\gamma) = \lambda$ is also sufficient for strong convergence, we argue as follows. Recall that any minimizing sequence $(\gamma_n,\alpha_n)$ with $\lambda < \lambda^{\rm HFB}(\kappa)$ is bounded in $\cX$, by Lemma \ref{lem:coercive}. Then Lemma \ref{lem:estim_density} implies that $\sqrt{\rho}_{\gamma_n}$ is bounded in $\Hhalf$. Thus we may assume that $\sqrt{\rho}_{\gamma_n} \wto \sqrt{\rho}_\gamma$ weakly in $H^{1/2}(\R^3)$ and $\rho_{\gamma_n} \to \rho_{\gamma}$ pointwise almost everywhere. Using this pointwise convergence and that $\int_{\R^3} \rho_{\gamma_n} < C$, we deduce from the Brezis-Lieb refinement of Fatou's lemma (see, e.\,g., \cite[Theorem 1.9]{LieLos-01}) that
$$
\int_{\R^3} \rho_{\gamma_n}(x) \, dx = \int_{\R^3} \rho_\gamma(x) \, dx + \int_{\R^3} |\rho_{\gamma}(x) - \rho_{\gamma_n}(x)| \, dx + o(1),
$$ 
where $o(1) \rightarrow 0$ as $n \rightarrow \infty$. Therefore $\tr(\gamma) = \int \rho_\gamma = \lambda$ implies that $\rho_{\gamma_n} \rightarrow \rho_\gamma$ strongly in $L^1(\R^3)$. Furthermore, by interpolation and the fact that $\| \rho_{\gamma_n} \|_{L^p} < C$ for $1 \leq p \leq 3/2$, we deduce that 
\begin{equation}
\mbox{$\rho_{\gamma_n} \to \rho_\gamma$ strongly in $L^p(\R^3)$ for $1 \leq p < 3/2$.}
\end{equation}
Thus the Hardy-Littlewood-Sobolev inequality shows that $D(\rho_{\gamma_n}, \rho_{\gamma_n}) \rightarrow D(\rho_\gamma, \rho_\gamma)$. The exchange term is found to be weakly-$*$ lower semicontinuous by Fatou's lemma:
\begin{equation}
\liminf_{n \rightarrow \infty} {\rm Ex}(\gamma_n) \geq {\rm Ex}(\gamma).
\end{equation} 
By Proposition \ref{prop:pairing_kinetic}, we have that the kinetic energy plus the pairing term form a wlsc-$\ast$ functional on $\cK_\lambda$. Hence, we conclude
\begin{equation}
\liminf_{n \rightarrow \infty} \cE(\gamma_n) \geq \cE(\gamma) \geq I(\lambda) .
\end{equation}
As by assumption $\cE(\gamma_n) \rightarrow I(\lambda)$, this proves {\em a posteriori} that $\gamma$ is a minimizer for $I(\lambda)$ and that
\begin{equation}
\lim_{n\to\ii}\cG(\gamma_n,\alpha_n)=\cG(\gamma,\alpha).
\end{equation}
By Proposition \ref{prop:pairing_kinetic} again, this shows that  $(\gamma_n,\alpha_n) \rightarrow (\gamma,\alpha)$ strongly in $\cX$. 
\end{proof}

Next, we study the minimization of the functional $\cG$ in more detail. This will be useful for the study of a minimizing sequence for $I(\lambda)$ when ruling out \emph{vanishing}; see Section \ref{sec:proof_thm_relaxed}. We start by introducing the following variational problem:
\begin{equation}
\boxed{
G(\lambda):=\inf\{\cG(\gamma,\alpha) : (\gamma,\alpha)\in\cK_\lambda\}.
}
\label{def_K}
\end{equation} 
Indeed, we can derive the following formula for this infimum.
\begin{prop} {\bf (Value of $G(\lambda)$ and upper bound on $I(\lambda)$)}. \label{prop:comparison_I_J}
Let $m>0$ and $0<\kappa<4/\pi$. Then for $0\leq\lambda\leq \lambda^{\rm HFB}(\kappa)$, we have the equality
\begin{equation}
G(\lambda)=\beta\lambda ,
\label{value_J}
\end{equation}
where $m(1-\kappa\pi/4) < \beta < m$ is given by 
$$
\beta = \left \{ \begin{array}{ll} \displaystyle \inf_{\psi \in L^2_{\mathrm{odd}}(\R^3), \| \psi \|_{L^2} = 1} \langle \psi, H \psi \rangle & \mbox{if $q =1$}, \\[1ex]
\displaystyle \inf_{\psi \in L^2(\R^3), \| \psi \|_{L^2} = 1} \langle \psi, H \psi \rangle & \mbox{if $q \geq 2$.} \end{array} \right .
$$
Here the integer $q \geq 1$ denotes the spin degrees of freedom, and $H = \sqrt{-\Delta + m^2} - \frac{\kappa}{2 |x|}$ is defined in the sense of quadratic forms on $L^2(\R^3)$ with form domain $H^{1/2}(\R^3)$.

Moreover, we have the strict inequality, for all $0 < \lambda < \lambda^{\rm HFB}(\kappa)$,
\begin{equation}
 I(\lambda)<G(\lambda)-m\lambda.
\label{estim_I_J}
\end{equation}
\end{prop}

The proof of Proposition \ref{prop:comparison_I_J} is given in Appendix \ref{appendix:F}. It is partly inspired of the discussion of the M\"uller energy functional in \cite{FraLieSeiSie-07}.

\begin{remark}\label{rmk:unbounded}
Adapting the proof of Proposition \ref{prop:comparison_I_J} given in Appendix \ref{appendix:F}, one can show that $G(\lambda)=I(\lambda)=-\ii$ if $\kappa>4/\pi$. Indeed, in this case the one-body operator $H$ is not bounded below anymore; see \cite{Herbst-77}.
\end{remark}

\section{Proof of Theorem \ref{thm:prop_min}} \label{sec:proof_prop_min}

In this section we prove Theorem \ref{thm:prop_min} which establishes key properties of any minimizer $(\gamma,\alpha)$ of the HFB energy functional $\cE$. In fact, we shall need some of these properties (more precisely, a certain decay estimate) in the proof by contradiction of Theorem \ref{thm:relaxed} that establishes the existence of minimizers; see Section \ref{sec:proof_thm_relaxed} below.

\subsection{Mean-Field Equations and Bounds on the Lagrange Multiplier $\mu$} 
Let $(\gamma,\alpha)\in\cK$ a minimizer for $I(\lambda)$. To derive the equation satisfied by $(\gamma, \alpha)$, we can follow a general argument given in \cite{BacLieSol-94}. That is, by writing that $\cE\big((1-t)(\gamma,\alpha)+t(\gamma',\alpha')\big)\geq \cE(\gamma,\alpha)$ for all $(\gamma',\alpha')\in\cK_\lambda$, we see that $(\gamma,\alpha)$ is a solution of the following {\em linearized variational problem}
\begin{equation}
K(\lambda)=\frac12\inf_{\substack{(\gamma',\alpha')\in\cK\\ \tr(\gamma')=\lambda}}\tr_{P^0_-} \left\{F_\Gamma\left(\Gamma'-P^0_-\right)\right\}  ,
\label{min_pb_GAMMA}
\end{equation}
where $\tr_{P^0_-} ( \cdot )$ denotes the $P^0_-$-{\em trace} as defined in \cite{HaiLewSer-05a}. We use the definition of $F_\Gamma$ and $\Gamma$ introduced in \eqref{def_H_SCF2} and \eqref{def_GAMMA}, respectively, and we employ the notation
\begin{equation}
P^0_-:=\left(\begin{matrix} 0&0\\ 0&1\end{matrix}\right), \quad \Gamma' := \left ( \begin{matrix} \gamma' & \alpha ' \\ (\alpha')^* & 1-\overline{\gamma'} \end{matrix} \right ).
\end{equation}
Note that
\begin{equation}
\tr(\gamma')=\frac12 \str \left\{N(\Gamma'-P^0_-)\right\} ,
\end{equation}
and that the function $\lambda\mapsto K(\lambda)$ defined in \eqref{min_pb_GAMMA} is easily seen to be convex. Hence $K(\lambda)$ satisfies
\begin{equation} \label{ineq:Klambda}
K(\lambda') \geq K(\lambda) + \mu (\lambda' - \lambda) , \quad \mbox{for all $\lambda' > 0$},
\end{equation}
whenever $\mu \in [K'_-(\lambda),K'_+(\lambda)]$, where $K'_-(\lambda)$ and $K'_+(\lambda)$ denote the left and right derivatives of $K(\lambda)$, respectively. (Note that $K'_\pm(\lambda)$ indeed exist, since $\lambda\mapsto K(\lambda)$ is convex.) Following the proof of Lemma \ref{lem:IHF}, we find that $K(\lambda)$ is non-increasing,  and thus $\mu\leq0$ holds. Moreover, by using \eqref{ineq:Klambda}, we deduce that $(\gamma,\alpha)$ is also a minimizer for the following problem without constraint on the trace:
\begin{equation}
\inf_{(\gamma',\alpha')\in\cK} \str \left\{(F_\Gamma-\mu N)\left(\Gamma'-P^0_-\right)\right\}  ,
\label{min_pb_GAMMA2}
\end{equation}
where $N$ is defined in \eqref{def_N}. This observation allows us to deduce (like in \cite{HaiLewSer-05a}) that $\Gamma=\Gamma(\gamma, \alpha)$ solves the equation
\begin{equation}
\Gamma=\chi_{(-\ii,0)}\left(F_\Gamma-\mu N\right)+D ,
\end{equation}
where $D$ is a finite rank operator of the same matrix form as $\Gamma$ and satisfies ${\rm ran} (D)\subset \ker (F_\Gamma-\mu N)$.
As an important next step in our analysis, we now prove an upper bound on the chemical potential $\mu$ . 

\begin{lemma}{\bf (Upper bound on $\mu$).} \label{lem:estim_multiplier}
Let $\Gamma = \Gamma(\gamma, \alpha)$ be a minimizer for $I(\lambda)$ and let $\mu$ be as above. Then we always have that $\mu<0$ holds. If moreover $\alpha \neq 0$, then $\mu < \beta-m$ holds, where $\beta < m$ is given by Proposition \ref{prop:comparison_I_J} above.
\end{lemma}

\begin{proof}
The proof in the case $\alpha=0$ (Hartree-Fock case) is standard and deserves no further discussion. Hence we will assume that $\alpha\neq0$ from now on.
We start by expanding the energy around the minimizer $\Gamma$ up to second order. We make use of Bogoliubov transformations of the form $e^{i\epsilon H}$ with
$$H=\left(\begin{matrix}
h&k\\ k^*&-\overline{h}
\end{matrix}\right),\qquad h^*=h,\quad k^T=-k, $$
where $h \in \gS_1$ and $k \in \gS_2$. From \cite{BacLieSol-94} we recall that, for every $\epsilon \in \R$, we have that $\Gamma_\epsilon=e^{i\epsilon H}\Gamma e^{-i\epsilon H}$ is an HFB state. Expanding up to second order, we find 
$$\Gamma_\epsilon=\Gamma+\epsilon\Gamma_1+\epsilon^2\Gamma_2+O(\epsilon^3),$$
where
$$\Gamma_1=i[H,\Gamma],\qquad \Gamma_2=H\Gamma H-\frac12 H^2\Gamma-\frac12 \Gamma H^2.$$
Next, we take $h=0$ and $k=-ia$ with $a=-a^T \in \gS_2$ to be chosen below, that is:
$$H=\left(\begin{matrix}
0&-i a\\ ia^*&0
\end{matrix}\right).$$
An elementary calculation yields
$$\Gamma_1=\left(\begin{matrix}
0&a\\ a^*&0
\end{matrix}\right)
+\left(\begin{matrix}
\gamma_1& \alpha_1\\
\alpha_1^*&-\overline{\gamma_1}
\end{matrix}\right),\qquad \Gamma_2=\left(\begin{matrix}
aa^*&0\\ 0&-a^*a
\end{matrix}\right)
+\left(\begin{matrix}
\gamma_2& \alpha_2\\
\alpha_2^*&-\overline{\gamma_2}
\end{matrix}\right)$$
with 
$$\gamma_1=a\alpha^*+\alpha a^*,\quad \alpha_1= -\gamma a-a\overline{\gamma},$$
$$\gamma_2=-a\overline{\gamma}a^*-\frac12aa^*\gamma-\frac12\gamma aa^*,\qquad 
\alpha_2=-a\alpha^* a-\frac12aa^*\alpha-\frac12\alpha a^*a.$$
Calculating the energy and using that $[F_\Gamma-\mu N,\Gamma]=0$, we find
$$\cE(\Gamma_\epsilon)=\cE(\Gamma)+\epsilon\mu\tr ( \gamma_1 ) +\epsilon^2\cQ(a,a)+O(\epsilon^3),$$
where
\begin{multline*}
\cQ(a,a)=\tr(T(a^*a+\gamma_2))-\kappa D(\rho_\gamma,\rho_{a^*a+\gamma_2})+ \kappa \Re X(\gamma,a^*a+\gamma_2)-\kappa \Re X(\alpha,\alpha_2)\\
-\frac\kappa2 D(\rho_{\gamma_1},\rho_{\gamma_1})+\frac\kappa2 X(\gamma_1,\gamma_1)-\frac\kappa2 X(a+\alpha_1,a+\alpha_1).
\end{multline*}
Here we use the convenient notation
\begin{equation}
X(a,b) = \iint_{\R^3 \times \R^3} \frac{\overline{a(x,y)} b(x,y)}{|x-y|} \, dx \, dy .
\label{eq:defX}
\end{equation}

We now claim the following fact.

\begin{lemma}\label{lem:exists_a}
There exists $a=-a^T$ with $\tr(Ka^*a)<\ii$ such that $\tr(\gamma_1)<0$ and $\cQ(a,a)<(\beta-m)\tr(\gamma_2)$.
\end{lemma}

Assuming Lemma \ref{lem:exists_a} for the moment, we have $\tr(\gamma_\epsilon)<\tr(\gamma)=\lambda$ for $\epsilon > 0$ small. Furthermore, with the help of Proposition \ref{prop:comparison_I_J}, we deduce
\begin{align*}
\cE(\Gamma_\epsilon)\geq I(\tr(\gamma_\epsilon))&\geq I(\lambda)-I\big(\lambda-\tr(\gamma_\epsilon)\big) \geq I(\lambda)-(\beta-m)\big(\lambda-\tr(\gamma_\epsilon)\big)\\
&= I(\lambda)+(\beta-m)\big(\epsilon\tr(\gamma_1)+\epsilon^2\tr(\gamma_2)+O(\epsilon^3)\big).
\end{align*}
Collecting these estimates, we obtain
$$(\beta-m-\mu)\tr(\gamma_1)\leq \epsilon\big(\cQ(a,a)-(\beta-m)\tr(\gamma_2)\big)+O(\epsilon^2),$$
whence the  desired estimate on $\mu$ follows, by taking $\epsilon > 0$ sufficiently small. This concludes the proof of Lemma \ref{lem:estim_multiplier}, provided that Lemma \ref{lem:exists_a} holds, which we will show next.
\end{proof}

\begin{proof}[Proof of Lemma \ref{lem:exists_a}] 
We begin with a preliminary observation: Assume that we have constructed an $a$ with $a^T=-a$ and $\tr(Ka^*a)<\ii$, such that $\cQ(a,a)-(\beta-m)\tr(\gamma_2)<0$. If $\tr(\gamma_1)>0$, it suffices to replace $a$ by $-a$ to conclude the proof. If $\tr(\gamma_1)=0$, then we can replace $a$ by $a'=a-\eta \alpha$. We then have (with an obvious notation) $\tr(\gamma_1')=-2\eta\tr(\alpha\alpha^*)<0$ and $\cQ(a',a')-(\beta-m)\tr(\gamma'_2)= \cQ(a,a)-(\beta-m)\tr(\gamma_2)+O(\eta)$. Hence we can make the last term $<0$, by choosing $\eta > 0$ sufficiently small. As a conclusion, we only have to construct $a$ such that $\cQ(a,a)<(\beta-m)\tr(\gamma_2)$ without considering the constraint $\tr(\gamma_1)<0$.

To construct $a$, we use trial states in the spirit of Appendix \ref{appendix:F} below. Here, let us for simplicity assume that we have $q=1$ spin degrees of freedom (since the proof for $q\geq 2$ is analogous). That is, we choose 
$$a_{R,L}(x,y)=\chi_L(R(x-3L\vec{v}))f(Rx-Ry)\chi_L(R(y-3L\vec{v}))$$ 
where $R\in SO(3)$ is a rotation and $\vec{v} \in \R^3$ is a fixed normalized vector; and $f \in H^{1/2}(\R^3)$ is the normalized ground state solution to 
$(K-\kappa/2|x|)f=\beta f$,
with $K - \kappa/2|x|$ acting on $L^2_\mathrm{odd}(\R^3)$. Moreover, for $L > 0$ given, we define $\chi_L(x)=L^{-3/4}\chi(x)$ with $\chi \in C^\infty_0(\R^3)$ being a radial function with support in the unit ball $\{ |x| \leq 1 \}$. 
We will use the following fact:
\begin{equation}
\norm{K^{1/2}a_{R,L}}\leq \frac{C}{L^{3/2}},
\label{estim:K_a_L}
\end{equation}
where $\| \cdot \|$ denotes the operator norm on $L^2$. Assuming this bound, it is straightforward to derive that
$$\norm{K^{1/2}\gamma_1}_{\gS_2}+ \norm{K^{1/2}\alpha_1}_{\gS_1}\leq\frac{C}{L^{3/2}},$$
$$\tr(K\gamma_2)+  \norm{K^{1/2}\gamma_2}_{\gS_1}+ \norm{K^{1/2}\alpha_2}_{\gS_2}\leq\frac{C}{L^3}.$$
Using these bounds and the fact that $|X(f,g)|\leq C\norm{K^{1/2}f}_{\gS_2}\norm{K^{1/2}g}_{\gS_2}$, we deduce that
\begin{multline*}
\cQ(a_{R,L},a_{R,L})-(\beta-m)\tr(\gamma_2)\leq \tr\big((K-\beta)a_{R,L}a_{R,L}^*\big)-\frac\kappa2 X(a_{R,L},a_{R,L})\\
-\frac\kappa2 D(\rho_\gamma,\rho_{a_{R,L}a_{R,L}^*})+O\left(L^{-3/2}\right). 
\end{multline*}
The first term of the right-hand side is invariant under translation and rotations (and thus independent of $R \in SO(3)$). The estimates in Appendix \ref{appendix:F} give
$$\tr\big((K-\beta)a_{R,L}a_{R,L}^*\big)-\frac\kappa2 X(a_{R,L},a_{R,L})=O\left(L^{-2}\right),$$
which implies
$$\cQ(a_{R,L},a_{R,L})-(\beta-m)\tr(\gamma_2)\leq -\frac\kappa2 D(\rho_\gamma,\rho_{a_{R,L}a_{R,L}^*})+O\left(L^{-3/2} \right). $$
Now we average over $R\in SO(3)$ and apply Newton's Theorem to obtain
$$\int_{SO(3)}dR\;D(\rho_\gamma,\rho_{a_{R,L}a_{R,L}^*})\geq \frac{\tr(a_{1,L}a_{1,L}^*)\int_{|x|\leq L}\rho_\gamma}{4L}\geq\frac\epsilon{L},$$
for $L > 0$ large enough and some constant $\epsilon >0$. Hence the conclusion of Lemma \ref{lem:exists_a} follows for a suitable rotation $R=R(L)$.

It remains to prove the key estimate \eqref{estim:K_a_L}. First, by rotational and translation symmetry, we can assume that $R=1 \in SO(3)$ and $v =0 \in \R^3$ in the definition of $a_{R,L}$ above. Thus it suffices to consider $
a_L(x,y) = \chi_L(x) f(x-y) \chi_L(y)$ and to show that
$\| K^{1/2} a_L \| \leq {C}/{L^{3/2}}$.
To prove this, we note that $a_L$ acts like $a_L \psi = \chi_L (f  \ast (\chi_L \psi))$, for $\psi \in L^2(\R^3)$.
Next, we observe that
\begin{equation}
K^{1/2} a_L \psi = [K^{1/2}, \chi_L] f \ast (\chi_L \psi) + \chi_L ((K^{1/2} f )\ast ( \chi_L \psi ))
\label{eq:aL}
\end{equation}
Since $f(x)$ decays exponentially (which, e.\,g., follows from the method in \cite{DalOstSto-08}), we conclude
\begin{equation}
\| f \|_{L^1} \leq C, \qquad \| K^{1/2} f \|_{L^1} \leq C.
\label{ineq:Kf}
\end{equation}
Indeed, the first bound follows immediately from the exponential decay of $f(x)$. To derive the second bound, we rewrite the equation satisfied by $f$ to find that
$$
K^{1/2} f = K^{-1/2} \left( \frac{\kappa}{2 |x|} f + \beta f \right)  .
$$
Note that the integral kernel of $K^{1/2} = (-\Delta + m^2)^{-1/4}$ with $m > 0$ belongs to $L^1(\R^3)$; see, e.\,g., \cite[p. 132]{Stein-70}. Since $\| |x|^{-1} f + \beta f\|_{L^1} < \infty$ (by exponential decay of $f$), the second bound in \eqref{ineq:Kf} follows.
Going back to \eqref{eq:aL}, we now obtain that
\begin{align*}
\| K^{1/2} a_L \psi \|_{L^2} & \leq  \| f \|_{L^1} \| [K^{1/2}, \chi_L] \|  \| \chi_L \|_{L^\infty} \| \psi \|_{L^2}   + \| K^{1/2} f \|_{L^1} \| \chi_L \|_{L^\infty}^2 \| \psi \|_{L^2} \\
& \leq C ( L^{-5/2} + L^{-3/2}) \| \psi \|_{L^2} .
\end{align*}
Here we used that $\| [K^{1/2}, \chi_L] \| \leq C \|\nabla \chi_L \|_{L^\infty} \leq C L^{-1-3/4}$ and $\| \chi_L \| = \| \chi_L \|_{L^\infty} \leq C L^{-3/4}$.  This implies the bound \eqref{estim:K_a_L} and completes the proof of Lemma \ref{lem:exists_a}. \end{proof}

\subsection{Decay estimate}
The next step consists in proving the following decay result for the density function of minimizers.

\begin{lemma} \label{lem:rho_decay}
For all $R > 0$ sufficiently large, we have
\begin{equation}
\tr ( T \zeta_R \gamma \zeta_R) + \int \zeta_R^2(x) \rho_\gamma(x) \, dx \leq \frac{C}{R^2}, 
\label{eq:decay_R}
\end{equation}
where $\zeta_R$ is the cutoff function defined in \eqref{eq:cutoff} and $C > 0$ denotes some constant independent of $R$. 
\end{lemma} 

\begin{proof}

When $\alpha\equiv0$, we have $\mu<0$ by Lemma \ref{lem:estim_multiplier}, hence $\gamma=\chi_{(-\ii,0)}(H_\gamma-\mu)+D$ is finite rank. Exponential decay is obtained by following the method in \cite{DalOstSto-08}. 

We now assume that $\alpha\neq0$. By Lemma \ref{lem:estim_multiplier}, we have $\mu<\beta-m$. Writing that $(F_\Gamma-\mu N)\Gamma\leq0$ and extracting the first diagonal term yields
\begin{equation}
A:=(H_\gamma-\mu)\gamma-\kappa X_\alpha \alpha^*\leq 0,
\label{ineq:gamma}
\end{equation}
where $X_\alpha$ denotes the operator with kernel $X_\alpha(x,y) = \alpha(x,y)/|x-y|$. 

On the other hand, we consider $[F_\Gamma-\mu N,\Gamma]=0$ and extract the upper right corner. This yields the following two-body equation for the pairing wavefunction $\alpha = \alpha(x,y) \in H^{1/2}(\R^3 \times \R^3)$:
\begin{equation}
\left((H_\gamma)_x+(H_\gamma)_y-\frac{\kappa}{|x-y|}-2\mu\right)\alpha= - \kappa(\gamma_x+\gamma_y)X_\alpha,
\label{eq:alpha} 
\end{equation}
Here $A_x$ means that $A$ as an operator on $L^2(\R^3)$ acts on the $x$-variable in $\alpha(x,y)$. Likewise $A_y$ means that $A$ acts on the $y$-variable in $\alpha(x,y)$.
Our goal is now to combine \eqref{ineq:gamma} and $\eqref{eq:alpha}$ in a useful way to extract decay information. 

\subsubsection*{{\bf Step 1: Estimates from \eqref{ineq:gamma}}}
Writing $\tr(\zeta_R(A+A^*)\zeta_R)\leq0$ where $A=A^*$ is defined in \eqref{ineq:gamma}, we deduce that
\begin{multline}
\tr\left(\frac{(\zeta_R)^2 T+T(\zeta_R)^2}{2} \gamma \right) - \mu \int_{\R^3} \zeta^2_R(x) \rho_\gamma(x) \, dx\\
\leq \kappa \iint_{\R^3 \times \R^3} \zeta^2_R(x) \frac{|\alpha(x,y)|^2}{|x-y|} \, dx \, dy 
+ \kappa \int V(x) \zeta^2_R(x) \rho_\gamma(x)  \, dx,
\end{multline}
where we put $V(x) = (|x|^{-1} \ast \rho_\gamma)(x)$. Note that we discarded the contribution of the exchange term in $H_\gamma$, because it gives a nonnegative contribution to the left-hand side of the above inequality. 

Since $\sqrt{\rho_\gamma} \in H^{1/2}(\R^3)$, one can see that $V(x) \to 0$ as $|x| \to \infty$. Hence 
 the term $\int V \zeta_R^2 \rho_\gamma$ can be estimated by $\delta(R) \int \zeta^2_R \rho_\gamma$ where $\delta(R)=\sup_{|x|\geq R}|V(x)| \to 0$ as $R \to \infty$. Moreover, the localization formula \eqref{eq:lemloc} from Lemma \ref{lem:IMS2} yields
\begin{multline}
\tr(T \zeta_R \gamma \zeta_R ) - \mu \int_{\R^3} \zeta^2_R(x) \rho_\gamma(x) \, dx \leq \delta(R) \int_{\R^3} \zeta^2_R(x) \rho_\gamma(x) \, dx \\  
+ \kappa \iint_{\R^3 \times \R^3} \zeta^2_R(x) \frac{|\alpha(x,y)|^2}{|x-y|} \, dx \, dy
 +\frac1\pi\tr\left( \left(\int_0^\ii \frac{1}{K^2+s}|\nabla\zeta_R|^2\frac{1}{K^2+s}\sqrt{s}\,ds\right)\gamma\right).
\label{eq:1st_estim_gamma}
\end{multline}
The last term can be estimated as usual by $\leq C\lambda/R^2$. 
Using that in particular $\mu < 0$, we finally deduce
\begin{equation}
\tr(T \zeta_R \gamma \zeta_R ) + \int_{\R^3} \zeta^2_R(x) \rho_\gamma(x) \, dx \leq  \frac{C}{R^2} 
+ C \iint_{\R^3 \times \R^3} \zeta^2_R(x) \frac{|\alpha(x,y)|^2}{|x-y|} \, dx \, dy ,
\label{ineq:master1}
\end{equation}
with $C > 0$ some constant and for $R>0$ sufficiently large. We will now estimate the right hand side of \eqref{ineq:master1} by using the two-body equation \eqref{eq:alpha}.

\subsubsection*{{\bf Step 2: Estimates from \eqref{eq:alpha}}}
We multiply equation \eqref{eq:alpha} by $\zeta_R(x) \zeta_R(y)$ from the left and we project onto
\begin{equation}
\alpha_R(x,y) = \zeta_R(x) \alpha(x,y) \zeta_R(y)  .
\end{equation}
This gives, where $\langle \cdot, \cdot \rangle$ is the inner product on $L^2(\R^3 \times \R^3)$, the following
\begin{multline*}
\langle \alpha_R, (T_x + T_y) \alpha_R \rangle - \kappa \iint_{\R^3 \times \R^3} \frac{|\alpha_R(x,y)|^2}{|x-y| } \, dx \, dy - 2 \mu \iint_{\R^3 \times \R^3} |\alpha_R(x,y)|^2 \, dx \, dy \\
=  I + II + III +IV,
\end{multline*}
where
$$
I= \langle \alpha_R, [ T_x + T_y , \zeta_R(x) \zeta_R(y)] \alpha \rangle,\quad
II = 2\kappa \iint V(x) \zeta^2_R(x) \zeta^2_R(y) |\alpha(x,y)|^2 \, dx \,dy ,$$
$$III = - \kappa \langle \alpha_R, \zeta_R(x) \zeta_R(y) ((X_\gamma)_x + (X_\gamma)_y) \alpha \rangle,\ \ 
IV = \kappa \langle \alpha_R, \zeta_R(x) \zeta_R(y) (\gamma_x + \gamma_y)) X_\alpha \rangle .$$
First, we note that $II$ is easy to estimate by
\begin{equation}
|II| \leq \delta(R) \int_{\R^3} \zeta^2_R(x) \rho_\gamma(x) \, dx,
\end{equation} 
with $\delta(R) \to 0$ as $R \to \infty$, by using that $\rho_{\alpha \alpha^*} \leq \rho_{\gamma}$ and $V \to 0$ as $|x| \to \infty$.

Next, we claim that
\begin{equation}
| III | + |IV|\leq \delta(R)\left ( \int_{\R^3} \zeta^2_R(x) \rho_\gamma(x) \, dx + \tr ( T \zeta_R \gamma \zeta_R ) \right ) .
\label{ineq:II}
\end{equation}
For $III$, this can be seen as follows (we only consider the term in $III$ that contains $(X_\gamma)_x$, the other being treated similarly):
\begin{multline}
| \langle \alpha_R, \zeta_R(x) \zeta_R(y) (X_\gamma)_x \alpha \rangle | =  | \tr( \alpha_R^*K^{1/2} K^{-1/2} X_{\zeta_R \gamma} \alpha \zeta_R )| \\
\leq \| K^{1/2} \alpha_R \|_{\gS_2} \| K^{-1/2} X_{\zeta_R \gamma} \|_{\gS_2} \| \alpha \zeta_R \|_{\gS_2},
\label{ineq:Xg}
\end{multline} 
using that $\| A \|_{\gS_\infty} \leq \| A \|_{\gS_2}$ holds. Next, by the Hardy-Kato inequality, we have that $\| K^{-1/2} X_f \|_{\gS_2} \leq C \| K^{1/2} f \|_{\gS_2}$. Therefore we find (using that $\alpha \alpha^* \leq \gamma$ and $\gamma^2 \leq \gamma$ and $\zeta_R^2 \leq 1$) the following bound:
\begin{align*}
|\eqref{ineq:Xg}| &\leq C \| K^{1/2} \alpha_R \|_{\gS_2} \| K^{1/2} \zeta_R \gamma \|_{\gS_2} \| \alpha \zeta_R \|_{\gS_2} \\
&\leq C \sqrt {\tr(K \zeta_R \alpha_R \zeta^2_R \alpha^*_R \zeta_R )} \sqrt{\tr(K \zeta_R \gamma^2 \zeta_R)} \sqrt{ \tr (K \zeta_R \alpha_R \alpha_R^* \zeta_R)} \\
&\leq C \tr(K \zeta_R \gamma \zeta_R)^{3/2} \leq \delta(R) \left ( \int_{\R^3} \zeta^2_R(x) \rho_\gamma(x) \, dx + \tr ( T \zeta_R \gamma \zeta_R ) \right ),
\end{align*}
where $\delta(R)\to0$ as $R\to\ii$. The proof is the same for $IV$.

Finally, we estimate $I$ as follows (we again only consider the term with $T_x$):
\begin{multline*}
|\langle \alpha_R, [T_x, \zeta_R(x)] \alpha \zeta_R(y) \rangle | \leq \| [T_x, \zeta_R ] \| \| \alpha_R \|_{\gS_2} \| \alpha \zeta_R \|_{\gS_2} \\ \leq \frac{C}{R} \tr(\zeta_R\gamma \zeta_R)  \leq \frac{C}{R} \int_{\R^3} \zeta^2_R(x) \rho_\gamma(x) \, dx,
\end{multline*}
where we used the commutator estimate $ \| [T_x, f(x)] \| \leq C \| \nabla f\|_{L^\infty}$ combined with the fact that $\alpha \alpha^* \leq \gamma$.

In summary, we have proved the following estimate
\begin{multline}
\langle \alpha_R, (T_x + T_y) \alpha_R \rangle - \kappa \iint_{\R^3 \times \R^3} \frac{|\alpha_R(x,y)|^2}{|x-y| } \, dx \, dy \\ - 2 \mu \iint_{\R^3 \times \R^3} |\alpha_R(x,y)|^2 \, dx \, dy 
\leq \delta(R) \left ( \int_{\R^3} \zeta^2_R(x) \rho_\gamma(x) \, dx + \tr ( T \zeta_R \gamma \zeta_R ) \right )  .
\end{multline}
Next, we turn to the left-hand side and derive a lower bound as follows. Recall from Appendix \ref{appendix:F} the lower bound
\begin{equation}
\left\langle \alpha, \left (T_x + T_y - \kappa \frac{1}{|x-y|} \right ) \alpha \right\rangle \geq 2 (\beta-m) \langle \alpha, \alpha \rangle
\end{equation}
for all $\alpha \in H^{1/2}(\R^3,\C^q)\wedge H^{1/2}(\R^3,\C^q)$, provided that $\kappa < 4/\pi$ holds. Note that $\beta = \beta(\kappa)$ depends continuously on $\kappa$.
Using now the strict inequality $2\mu < 2(\beta-m)$, we deduce that $2 \mu < 2\beta(\kappa+\epsilon) -2m$ for $\epsilon > 0$ small, and hence
\begin{multline}
\langle \alpha_R, (T_x + T_y) \alpha_R \rangle - (\kappa+\epsilon-\epsilon) \iint_{\R^3 \times \R^3} \frac{|\alpha_R(x,y)|^2}{|x-y| } \, dx \, dy \\ - 2 \mu \iint_{\R^3 \times \R^3} |\alpha_R(x,y)|^2 \, dx \, dy 
\geq \epsilon \iint_{\R^3 \times \R^3} \frac{|\alpha_R(x,y)|^2}{|x-y| } \, dx \, dy ,
\end{multline}
for some $\epsilon > 0$ sufficiently small. Therefore, we arrive at the following estimate:
\begin{multline}
 \iint_{\R^3 \times \R^3} \frac{\zeta^2_R(x) \zeta^2_R(y) |\alpha(x,y)|^2}{|x-y| } \, dx \, dy \\ \leq \delta(R) \left ( \int_{\R^3} \zeta^2_R(x) \rho_\gamma(x) \, dx + \tr ( T \zeta_R \gamma \zeta_R ) \right )  
 \label{ineq:master2}
\end{multline}
with $\delta(R) \to 0$ as $R \to \infty$.

\subsubsection*{{\bf Step 3: Combining \eqref{ineq:master1} and  \eqref{ineq:master2}}}
We define the sequence $\{ R_n \}_{n=1}^\infty$ of radii given by
$R_n = 4^n$, for $n \geq1$. From \eqref{ineq:master1} we conclude for all $n$ large enough that
\begin{multline}
\tr(T \zeta_{R_{n+1}} \gamma \zeta_{R_{n+1}} ) + \int_{\R^3} \zeta^2_{R_{n+1}}(x)^2 \rho_\gamma(x) \, dx\\
  \leq \frac{C}{(R_{n+1})^2}+ C \iint_{\R^3 \times \R^3} \zeta^2_{R_{n+1}}(x) \frac{|\alpha(x,y)|^2}{|x-y|} \, dx \, dy.
\end{multline}
Next, we let $\chi_{R_{n}} = \sqrt{1-\zeta_{R_n}^2}$ so that $\zeta_{R_n}^2 + \chi_{R_n}^2 \equiv 1$ holds. Therefore we can split and estimate the pairing term in the inequality above as follows:
\begin{multline*}
\iint_{\R^3 \times \R^3} \zeta^2_{R_{n+1}}(x)  \frac{|\alpha(x,y)|^2}{|x-y|} \, dx \, dy  = \iint_{\R^3 \times \R^3} \zeta^2_{R_{n+1}}(x) \chi^2_{R_n}(y) \frac{|\alpha(x,y)|^2}{|x-y|} \, dx \, dy \\
 + \iint_{\R^3 \times \R^3} \zeta^2_{R_{n+1}}(x) \zeta^2_{R_n}(y) \frac{|\alpha(x,y)|^2}{|x-y|} \,dx \, dy  =: I + II .
\end{multline*}
Since $R_n = \frac{1}{4} R_{n+1}$, we deduce from the support properties of $\zeta_{R_{n+1}}(x)$ and $\chi_{R_n}(y)$ that
\begin{equation*}
I \leq \frac{C}{R_{n}} \int_{\R^3} \zeta^2_{R_{n+1}}(x) \rho_\gamma(x) \, dx ,
\end{equation*}
using $\rho_{\alpha \alpha^*}(x) = \int_{\R^3} |\alpha(x,y)|^2 \, dy \leq \rho_\gamma(x)$ thanks to $\alpha \alpha^* \leq \gamma$. Further, we notice that $\zeta_{R_{n+1}} \leq \zeta_{R_n}$ and thus
\begin{equation*}
II \leq \iint \zeta^2_{R_n}(x) \zeta^2_{R_n}(y) \frac{|\alpha(x,y)|^2}{|x-y|} \, dx \, dy \leq \delta(R_n) \left ( \int \zeta^2_{R_n}(x) \rho_\gamma + \tr (T \zeta_{R_n} \gamma \zeta_{R_n}) \right ),
\end{equation*}
where we used inequality \eqref{ineq:master2} in the last step.

Let us now define the sequence $\{ u_n\}_{n=1}^\infty$ of nonnegative numbers given by
\begin{equation}
u_n = \int_{\R^3} \zeta_{R_n}^2(x) \rho_\gamma(x) \, dx + \tr(T \zeta_{R_n} \gamma \zeta_{R_n}) .
\end{equation}
Collecting the previous estimates, we obtain the recursive inequality:
\begin{equation}
u_{n+1} \leq \delta_n u_n + \frac{C}{(R_{n+1})^2} \quad \mbox{for $n \geq n_0$}.
\label{ineq:nn}
\end{equation}
where $\delta_n \to 0$ as $n \to \infty$ and $n_0 \geq 1$ is sufficiently large. A simple induction argument then shows that $u_n$ satisfies the following bound:
\begin{equation}
u_n \leq \frac{B}{(R_n)^2} \quad \mbox{for $n \geq 1$},
\label{ineq:un_bound}
\end{equation}
where $B > 0$ is some sufficiently large constant.

It remains to extend the bound in Lemma \ref{lem:rho_decay} to all $R>0$ sufficiently large. To this end, let $R\geq4$ and $n\geq1$ such that $4^n\leq R <4^{n+1}$.  From \eqref{ineq:master1} and the fact that $\zeta_R\leq \zeta_{4^{n-1}}$, we have
$$\tr(T \zeta_R \gamma \zeta_R ) + \int_{\R^3} \zeta_R^2(x) \rho_\gamma(x) \, dx \leq  \frac{C}{R^2} 
+ C \iint_{\R^3 \times \R^3} \zeta^2_{R_{n-1}}(x) \frac{|\alpha(x,y)|^2}{|x-y|} \, dx \, dy.$$
Our previous estimates show that 
$$\iint_{\R^3 \times \R^3} \zeta^2_{R_{n-1}}(x) \frac{|\alpha(x,y)|^2}{|x-y|} \, dx \, dy\leq \frac{C}{4^{2(n-1)}}\leq \frac{4^4C}{R^2},$$
and the desired result follows for $R>0$ sufficiently large. The proof of Lemma \ref{lem:rho_decay} is now complete. \end{proof}

\begin{remark}\label{rmk:better_decay}
One can bootstrap the decay estimates obtained above by commuting $|\nabla\zeta_R|^2$ with $(K^2+s)^{-1}$ on the right side of \eqref{eq:1st_estim_gamma}, leading to a bound $\leq C_k/R^k$ for all $k\geq1$ in \eqref{eq:decay_R}. This better decay estimate is, however, unnecessary for our existence proof, and hence we will not give any details here.
\end{remark}

Having established Lemma \ref{lem:rho_decay}, the proof of Theorem \ref{thm:prop_min} is now complete.  \hfill $\blacksquare$

\section{Proof of Theorem \ref{thm:infinite_rank}}\label{sec:proof_infinite_rank}
 
First, we recall the equality \eqref{eq:tau_min} and we note that, by assumption, we have that $\sqrt{\tau(1-\tau)} \neq 0$ holds. To show that $\tau$ must have an infinite rank, we note the similarity between $\mathcal{F}(\tau)$ and the M\"uller functional studied in \cite{FraLieSeiSie-07}. Indeed, by following an argument in \cite{FraLieSeiSie-07}, we can prove that $\tau$ has infinite rank as follows. For the reader's convenience, we provide the details of the adaptation.

We argue by contradiction and assume that $\tau$ is finite rank. Hence we can write
$$\tau=\sum_{i=1}^Kn_i|\phi_i\rangle\langle\phi_i|,\qquad 0< n_i\leq 1,\quad \sum_{i=1}^Kn_i=\lambda/2$$
with $\overline{\phi_i}=\phi_i$ normalized in $L^2(\R^3,\R)$. Moreover, we can assume $0<n_1<1$ without loss of generality since $\alpha \neq 0$ by assumption.
\begin{lemma}\label{lem:Sobolev}
We have $\phi_i\in H^s(\R^3,\R)$ for all $i=1..K$ and all $s\geq0$.
\end{lemma}
\begin{proof}[Proof of Lemma \ref{lem:Sobolev}]
Note that $\phi_1,...,\phi_K$ solve the following minimization problem
\begin{multline}
\inf\Bigg\{\sum_{i=1}^Nn_i\pscal{T\phi_i,\phi_i}-\kappa\iint\frac{\left(\sum_{i=1}^Nn_i\phi_i(x)^2\right)\left(\sum_{i=1}^Nn_i\phi_i(y)^2\right)}{|x-y|}dx\,dy\\+\frac\kappa2\iint\frac{\left(\sum_{i=1}^Nn_i\phi_i(x)\phi_i(y)\right)^2}{|x-y|}dx\,dy
-\frac\kappa2\iint\frac{\left(\sum_{i=1}^N\sqrt{n_i(1-n_i)}\phi_i(x)\phi_i(y)\right)^2}{|x-y|}dx\,dy,\\
(\phi_1,...,\phi_K)\in H^{1/2}(\R^3,\R),\ \pscal{\phi_i,\phi_j}=\delta_{ij}\Bigg\}.
\end{multline}
Hence they must solve the following nonlinear equation, 
\begin{multline}
n_i\Bigg(T-2\kappa\left(\sum_{k=1}^Nn_k\phi_k(x)^2\ast\frac{1}{|x|}\right)+\kappa\frac{\sum_{k=1}^Nn_k\phi_k(x)\phi_k(y)}{|x-y|}\\
-\kappa\sqrt{\frac{1-n_i}{n_i}}\frac{\sum_{k=1}^N\sqrt{n_k(1-n_k)}\phi_k(x)\phi_k(y)}{|x-y|}\Bigg)\phi_i=\sum_{k=1}^N\lambda_{ik}\phi_i
\end{multline}
where the $\lambda_{ij}$ are Lagrange multipliers associated with the constraints $\pscal{\phi_i,\phi_j}=\delta_{ij}$.
The result is then proved using $n_i>0$ and a simple bootstrap argument.
\end{proof}

Next, let $\epsilon>0$  be small enough and consider a real function $\psi\in{\rm span}(\phi_i,\ i=1..K)^\perp$. We introduce the following test state
$$\tau_\epsilon=(n_1-\epsilon)|\phi_1\rangle\langle\phi_1|+\sum_{i=2}^Kn_i|\phi_i\rangle\langle\phi_i|+\epsilon|\psi\rangle\langle\psi|$$
and compute
\begin{equation}
\cF(\tau_\epsilon)= \cF(\tau)-\kappa\sqrt{\epsilon} \iint_{\R^3\times\R^3}\frac{a(x,y)\psi(x)\psi(y)}{|x-y|}dx\,dy+O(\epsilon),
\end{equation}
where $a:=\sqrt{\tau(1-\tau)}\neq0$ is a (real) finite rank symmetric operator acting on $L^2(\R^3)$.
Hence we will get a contradiction once we can prove the following:
\begin{lemma}
There exists $\psi\in{\rm span}(\phi_i,\ i=1..K)^\perp$, $\overline{\psi}=\psi$, such that
$$\iint_{\R^3\times\R^3}\frac{a(x,y)\psi(x)\psi(y)}{|x-y|}dx\,dy=4\pi\sum_{i=1}^K\sqrt{n_i(1-n_i)}\int_{\R^3}\frac{\left|\widehat{\phi_i\psi}(k)\right|^2}{|k|^2}dk>0.$$
\end{lemma}
\begin{proof}
Assume on the contrary that
$$\iint_{\R^3\times\R^3}\frac{a(x,y)\psi(x)\psi(y)}{|x-y|}dx\,dy=0$$
for all $\psi\in {\rm span}(\phi_i,\ i=1..K)^\perp$. This means that the (real) symmetric Hilbert-Schmidt operator $R$ with kernel $a(x,y)|x-y|^{-1}$ vanishes when it is restricted to ${\rm span}(\phi_i,\ i=1..K)^\perp$. As $R\geq0$, this implies that we must have
\begin{equation}
 \frac{a(x,y)}{|x-y|}=\sum_{i,j=1}^Ka_{ij}\phi_i(x)\phi_j(y)
\label{equality_contradiction}
\end{equation}
where $(a_{ij})$ is a real non-negative symmetric matrix. Multipliying \eqref{equality_contradiction} by $|x-y|$ and taking $x=y$ we deduce that $a(x,x)=\sum_{i=1}^K\sqrt{n_i(1-n_i)}\phi_i(x)^2=0$ for every $x\in\R^3$. This implies $a=\sqrt{\tau(1-\tau)}=0$ which is a contradiction with $\tau^2\neq\tau$.
\end{proof}

The proof of Theorem \ref{thm:infinite_rank} is now complete.  \hfill $\blacksquare$
 
\section{Proof of Theorem \ref{thm:relaxed}}\label{sec:proof_thm_relaxed}
In this section, we present the proof of Theorem \ref{thm:relaxed}. By Lemma \ref{lem:IHF}, we note that Part (ii) will be an immediate consequence of Part (i). Hence it suffices to prove Part (i), where we argue by contradiction as follows. We suppose that
$0 < \lambda < \lambda^{\rm HFB}(\kappa)$,
and we assume throughout the proof that there exists a minimizing sequence $\{(\gamma_n,\alpha_n)\}_{n \in \N}\subset\cK_\lambda$ for $I(\lambda)$ that is not relatively compact in $\cX$ up to translations. We divide our proof into several subsections as follows.

\subsubsection*{\bf Step 1: No Vanishing}
We begin by ruling out vanishing of the sequence $\{ \gamma_n \}_{n \in \N}$. By vanishing, we mean that the associated sequence of densities $\{ \rho_{\gamma_n} \}_{n \in \N}$ satisfies
\begin{equation} \label{eq:vanish}
\limsup_{n \rightarrow \infty} \left ( \sup_{y \in \R^3} \int_{|x-y|\leq R}\rho_{\gamma_n}(x)\, dx\right) = 0,
\end{equation}
for all $R >0 $. 
For each density matrix $\gamma_n$, we define the following concentration function \cite{Lions-84}
\begin{equation} 
Q^1_n(R)= \sup_{y \in \R^3} \int_{|x-y|\leq R}\rho_{\gamma_n}(x)\, dx,
\label{eq:Q1n}
\end{equation}
Furthermore, we introduce the number
\begin{equation} 
\lambda^1:=\lim_{R \to \infty} \limsup_{n \to \infty} Q^1_n(R),
\end{equation}
Indeed, we shall rule out that $\lambda_1 = 0$ holds. This implies that  the sequence $\{ \gamma_n \}$ cannot vanish in the sense defined above.

\begin{lemma}{\bf (No Vanishing).} \label{lem:novanish} 
We have $\lambda^1>0$.
\end{lemma}

To prepare the proof of Lemma \ref{lem:novanish}, we first establish the following fact, which is an adaptation of a similar result in \cite[Lemma I.1]{Lions-84} to the fractional Sobolev space $H^{1/2}(\R^3)$. An essential ingredient will be the localization formula in Lemma \ref{lem:IMS1}.

\begin{lemma} \label{lem:vanish}
Let $\{ f_n \}$ be a bounded sequence in $\Hhalf$ and suppose $\{ |f_n|^2 \}$ vanishes in the sense that
\[
\limsup_{n \rightarrow \infty} \left ( \sup_{y \in \R^3}  \int_{|x-y| \leq R} |f_n(x)|^2 \, dx  \right ) = 0, \quad \mbox{for all $R >0$.}
\]
Then $\| f_n \|_{L^p} \to 0$ for $2 < p < 3$.
\end{lemma}

\begin{proof}
We consider a smooth partition of unity $\{ \chi_k \}_{k \in \Z^3}$ with $\sum_{k \in \Z^3} \chi_k^2 \equiv 1$ such that
$\chi_k \equiv 1$ on $C_k$ and $\mathrm{supp} \, \chi_k \subset C'_k$, 
with the half-open cubes $C_k = [k,k+1)^3$ and $C'_k = [k-1,k+2)^3$ where $k \in \Z^3$. Furthermore, we assume that $\sup_{k \in \Z^3} |\nabla \chi_k(x)| \leq C$. Next, we estimate as follows.
\begin{align}
\int_{\R^3} |f_{n}(x)|^{8/3} \, dx & = \sum_{k \in \Z^3} \int_{\R^3} \chi_k(x)^2 |f_{n}(x)|^{8/3} \, dx  \nonumber \\
& = \sum_{k \in \Z^3} \int_{\R^3} |\1_{C'_k}(x)  f_{n}(x)|^{2/3} | \chi_k(x) f_{n}(x)|^2 \, dx \nonumber \\
& \leq \sum_{k \in \Z^3} \big \| \1_{C'_k} |f_n|^2 \big \|_{L^1(\R^3)}^{1/3} \big \| \chi_k f_{n} \big \|_{L^3(\R^3)}^2 \nonumber \\
& \leq C \sup_{k \in \Z^3} \big \| \1_{C'_k} |f_n|^2 \big \|_{L^1(\R^3)}^{1/3} \sum_{k \in \Z^3} \langle f_{n},  \chi_k K \chi_k f_{n} \rangle. \label{ineq:vanish}
\end{align}
We have used the Sobolev-type inequality $\| f \|_{L^3(\R^3)}^2 \leq C \langle f, K f \rangle$. By Lemma \ref{lem:IMS1}, the nonlocal operator $K$ satisfies
\begin{equation}
\sum_{k \in \Z^3} \chi_k K \chi_k   \leq K + \frac{1}{\pi} \int_0^\infty \frac{1}{s+K^2} \left ( \sum_{k \in \Z^3} |\nabla \chi_k |^2 \right ) \frac{1}{s+K^2} \, \sqrt{s} \, ds 
\leq K + C,
\end{equation}
since the integral expression on the right side in the first inequality is a bounded self-adjoint operator due to the fact that $\sum_{k \in \Z^3} | \nabla \chi_k(x)|^2 \in L^\infty(\R^3)$ by our choice of the partition $\{ \chi_k \}_{k \in \Z^3}$. Therefore, we conclude
\begin{equation}
\sum_{k \in \Z^3} \langle f_{n},  \chi_k K \chi_k f_{n} \rangle \leq C \| f_{n} \|_{\Hhalf}^2 .
\end{equation} 
Since, by assumption, we have$\| f_{n} \|_{\Hhalf} \leq C$ independent of $n$, estimate (\ref{ineq:vanish}) leads to
\begin{equation}
\int_{\R^3} |f_{n}(x)|^{8/3} \,dx \leq C \sup_{k \in \Z^3} \big \| \1_{C'_k} |f_{n}|^2 \big \|_{L^1(\R^3)}^{1/3}  \rightarrow 0,
\end{equation}
using that $\{ |f_n|^2 \}_{n \in \N}$ vanishes. This shows that $f_{n} \rightarrow 0$ in $L^{8/3}(\R^3)$, whence it converges to 0 strongly in $L^p(\R^3)$ for all $2 < p < 3$, by interpolation and the fact that $\{ f_{n} \}_{n \in \N}$ is a bounded sequence in $L^2(\R^3) \cap L^{3}(\R^3)$ due to Sobolev embedding. This completes the proof of Lemma \ref{lem:vanish}. \end{proof}

We are now ready to prove Lemma \ref{lem:novanish}.

\begin{proof}[Proof of Lemma \ref{lem:novanish}]
By the coercivity of $\cE$ stated in Lemma \ref{lem:coercive}, we know that our minimizing sequence $(\gamma_n,\alpha_n)$ is bounded in $\cX$. Hence we deduce from Lemma \ref{lem:estim_density} that $\sqrt{\rho}_{\gamma_n}$ is bounded in $H^{1/2}(\R^3)$. 
Applying Lemma \ref{lem:vanish} to the sequence $f_n = \sqrt{\rho}_{{\gamma}_n}$, we obtain that $\rho_{\gamma_n} \rightarrow 0$ strongly in $L^p(\R^3)$ for $1 < p < 3/2$. By the Hardy-Littlewood-Sobolev inequality and (\ref{estim_exchange}), this implies
\begin{equation}
D(\rho_{\gamma_n}, \rho_{\gamma_n}) \rightarrow 0 \quad \mbox{and} \quad {\rm Ex}(\gamma_n) \rightarrow 0.
\end{equation}
Therefore, we conclude
\begin{equation}
\cE(\gamma_n,\alpha_n)=\cG(\gamma_n,\alpha_n)-m\lambda+o(1)\geq G(\lambda)-m\lambda+o(1)
\label{estim_below_if_vanishing} 
\end{equation}
where we recall that $\cG$ and $G$ were defined in Section \ref{sec:auxillary}. But we have $\cE(\gamma_n,\alpha_n)\to I(\lambda)$ and $I(\lambda)<G(\lambda)-m\lambda$, by Proposition \ref{prop:comparison_I_J}. Thus \eqref{estim_below_if_vanishing} yields a contradiction and we must have $\lambda^1>0$.
\end{proof}

\subsubsection*{\bf Step 2: Dichotomy and Local Compactness}
We already know that our minimizing sequence $\{(\gamma_n,\alpha_n) \}_{n \in \N}$ with $\tr(\gamma_n) = \lambda$ does not vanish and we have assumed that it is not relatively compact in $\cX$ up to translations. By Corollary \ref{lem:compact}, this means that we must have $0<\lambda^1<\lambda$.
We will now use the following adaptation of a the classical dichotomy result \cite{Lions-84} to our operator setting. The precise statement reads as follows.

\begin{lemma}{\bf (Strong Local Convergence).} \label{lem:local_compactness} Up to extraction of a subsequence, there exists a density matrix $\gamma^1 \in \cK$ with $\tr(\gamma^1)=\lambda^1$, sequences $\{ R_n \}_{n \in \N} \subset \R_+$ with $R_n\to\ii$ and $\{y_n\}_{n \in\N} \subset \R^3$, such that
\begin{equation}
\sqrt{K} \tau_{y_n}^* \chi (R_n^{-1} \cdot )\gamma_n\chi (R_n^{-1} \cdot ) \tau_{y_n} \sqrt{K}\wto \sqrt{K}\gamma^1\sqrt{K}
\end{equation}
weakly-$\ast$ in $\gS_1$, where $K=\sqrt{-\Delta+m^2}$, and
\begin{equation}
\tau_{y_n}^* \chi (R_n^{-1} \cdot )\gamma_n\chi (R_n^{-1} \cdot ) \tau_{y_n} \to \gamma^1
\end{equation}
strongly in $\gS_1$. Here $\tau_y$ is unitary action representing translations in $\R^3$ as defined in \eqref{eq:tau}, and $0  \leq \chi \leq 1$ is the smooth cutoff function introduced in \eqref{eq:cutoff}. Moreover, we have
\begin{equation}
\lim_{n\to\ii} \int_{R_n\leq|x-y_n|\leq 2R_n}\rho_{\sqrt{K}\gamma_n\sqrt{K}}(x)\,dx=\lim_{n\to\ii} \int_{R_n\leq|x-y_n|\leq 2R_n}\rho_{\gamma_n}(x)\,dx=0.
\label{estim_density_conc_compt} 
\end{equation}
\end{lemma}
\begin{proof}
We will not detail the proof which is an adaptation of ideas by Lions \cite{Lions-84}, where one introduces another sequence of concentration functions
\begin{equation} 
Q^2_n(R)= \sup_{y \in \R^3} \int_{|x-y|\leq R}\rho_{\sqrt{K}\gamma_n\sqrt{K}}(x)\, dx 
\end{equation}
in addition to the sequence $\{ Q^1_n \}_{n=1}^\infty$ defined in \eqref{eq:Q1n}. We refer, for instance, to \cite{Friesecke-03} where a similar argument has been detailed. The rest of the proof follows well-known ideas of Lions \cite{Lions-84} coupled with the fact that our minimizing sequence $\{\gamma_n\}$ is bounded in $\cX$, hence (up to a subsequence) converges strongly locally in the trace class. The proof is the same in $\cX$.
\end{proof}

As we have already mentioned, our model given by $\cE(\gamma,\alpha)$ is invariant under spatial translations. Thus we may assume that (in Lemma \ref{lem:local_compactness}) the sequence of translations is given by
$$\boxed{y_n=0, \quad \mbox{for all $n \geq 1$}}$$ 
which we shall do for all the rest of the proof. 

\subsubsection*{\bf Step 3: Splitting of the Energy}
The next step shows that we must have 
\begin{equation} \label{eq:split}
I(\lambda)=I(\lambda^1)+I(\lambda-\lambda^1) \quad \text{and} \quad \mbox{$I(\lambda^1)$ has a minimizer.}
\end{equation}
To streamline our notation with Lemma \ref{lem:local_compactness}, we introduce a new radial cut-off function $0\leq\tilde\chi\leq1$ such that $\tilde{\chi} \equiv 1$ for $|x| \leq 7/5$ and $\tilde{\chi} \equiv 0$ for $|x| \geq 8/5$. Given the sequence $\{ R_n \}$ from Lemma \ref{lem:local_compactness}, we define the functions $\chi_n=\tilde\chi(\cdot/R_n)$ and $\zeta_{R_n}=\sqrt{1-\chi_n^2}$. Likewise, we define $\{(\gamma^1_n,\alpha^1_n) \}_{n \in \N}$ and $\{(\gamma^2_n,\alpha^2_n) \}_{n \in \N}$ by
$$(\gamma^1_n,\alpha^1_n) :=\chi_n(\gamma_n,\alpha_n) \chi_n \quad \mbox{and} \quad (\gamma^2_n,\alpha^2_n) :=\zeta_{R_n}(\gamma_n,\alpha_n)\zeta_{R_n}.$$
Note that, since $(\gamma_n, \alpha_n) \in \cK$, one easily verifies that $(\gamma^i_n,\alpha^i_n)\in\cK$ for $i=1,2$. Furthermore, we have $\gamma^1_n\to\gamma^1$ in $\gS_1$, by Lemma \ref{lem:local_compactness}, which in particular implies that $\tr(\gamma_n^1)=\int\rho_{\gamma_n}(x)\chi_n^2(x)\,dx\to\lambda^1$. 

Next, we shall prove that
\begin{equation} 
\cE(\gamma_n,\alpha_n) \geq \cE(\gamma_n^1,\alpha_n^1)+\cE(\gamma_n^2,\alpha_n^2)+o(1). 
\label{splitting}
\end{equation}
To deal with the kinetic energy, we use Lemma \ref{lem:IMS1} and $\tr (\gamma_n)=\lambda$ to find that
\begin{equation}
 \tr(T\gamma_n)\geq\tr(T\gamma^1_n)+\tr(T\gamma^2_n)-C\lambda R_n^{-2}.
\label{splitting_kinetic_energy}
\end{equation}
The next step in the proof of (\ref{splitting}) is to separate the direct term, which can be  done as follows:
\begin{equation*}
D(\rho_{\gamma_n},\rho_{\gamma_n})=D(\rho_{\gamma^1_n},\rho_{\gamma^1_n})+D(\rho_{\gamma^2_n},\rho_{\gamma^2_n})+2\iint\frac{\rho_{\gamma^1_n}(x)\rho_{\gamma^2_n}(y)}{|x-y|}dx\,dy.
\end{equation*}
Here, we note that 
\begin{align*}
&\iint\frac{\chi_n(x)^2\rho_{\gamma_n}(x)\zeta_n(y)^2\rho_{\gamma_n}(y)}{|x-y|}dx\,dy\\
&\qquad=\int_{|x|\leq 6R_n/5}\int_{|y|\geq 7R_n/5}\frac{\chi_n(x)^2\rho_{\gamma_n}(x)\zeta_n(y)^2\rho_{\gamma_n}(y)}{|x-y|}dx\,dy\\
&\qquad\qquad+\int_{6R_n/5\leq|x|\leq 8R_n/5}\int_{|y|\geq 9R_n/5}\frac{\chi_n(x)^2\rho_{\gamma_n}(x)\zeta_n(y)^2\rho_{\gamma_n}(y)}{|x-y|}dx\,dy\\
&\qquad\qquad+\int_{6R_n/5\leq|x|\leq 8R_n/5}\int_{7R_n/5\leq|y|\leq 9R_n/5}\frac{\chi_n(x)^2\rho_{\gamma_n}(x)\zeta_n(y)^2\rho_{\gamma_n}(y)}{|x-y|}dx\,dy .
\end{align*}
Therefore,
$$\iint\frac{\chi_n(x)^2\rho_{\gamma_n}(x)\zeta_n(y)^2\rho_{\gamma_n}(y)}{|x-y|}dx\,dy\leq \frac{10 \lambda^2}{R_n}+C\norm{\rho_{\gamma_n}}^2_{L^{6/5}(B(0,9R_n/5)\setminus B(0,6R_n/5))}.$$
Furthermore, by \eqref{estim_density_conc_compt}, we have that $\rho_{\gamma_n}\1_{B(0,2R_n)\setminus B(0,R_n)}\to0$ in $L^1(\R^3)$. Since this is also a bounded sequence in $L^p(B(0,2R_n)\setminus B(0,R_n))$ for all $1\leq p\leq3/2$, we infer by interpolation that it also tends to $0$ in $L^p(B(0,2R_n)\setminus B(0,R_n))$ for $1< p<3/2$. Hence, by the Hardy-Littlewood-Sobolev inequality, we conclude
$$D(\rho_{\gamma_n},\rho_{\gamma_n})=D(\rho_{\gamma^1_n},\rho_{\gamma^1_n})+D(\rho_{\gamma^2_n},\rho_{\gamma^2_n})+o(1).$$
As usual, the exchange term ${\rm Ex}(\gamma)$ is easily estimated by noticing that $ \Ex(\gamma_n)\geq\Ex(\gamma_n^1)+\Ex(\gamma_n^2)$.
The pairing term is treated similarly as done before several times
\begin{multline}
\iint\frac{\chi_n(x)^2\zeta_n(y)^2|\alpha_n(x,y)|^2}{|x-y|}dx\,dy\leq \frac{5}{R_n}\tr(\alpha_n^*\alpha_n)\\+
\int_{6R_n/5\leq|x|\leq 8R_n/5}\int_{7R_n/5\leq|y|\leq 9R_n/5}\frac{\chi_n(x)^2\zeta_n(y)^2|\alpha_n(x,y)|^2}{|x-y|}dx\,dy.
\label{last_ineg_pairing_localization} 
\end{multline}
Note that
\begin{equation}
\iint_{\R^3\times\R^3}|\alpha_n(x,y)|^2dx\,dy=\tr(\alpha_n\alpha^*_n)\leq \tr(\gamma_n)=\lambda,
\end{equation}
by \eqref{estim_pairing}. Define $\eta_n=\eta(\cdot/R_n)$ for some smooth function $\eta$ such that $\eta \equiv 1$ on the annulus $\{6/5\leq |x|\leq 9/5\}$ and $\eta \equiv 0$ outside the annulus $\{1\leq |x|\leq 2\}$. The last term of \eqref{last_ineg_pairing_localization} can be bounded above by
\begin{multline*}
\int_{6R_n/5\leq|x|\leq 8R_n/5}\int_{7R_n/5\leq|y|\leq 9R_n/5}\frac{\chi_n(x)^2\zeta_n(y)^2|\alpha_n(x,y)|^2}{|x-y|}dx\,dy\\
\leq \frac\pi2\tr\big(K\eta_n\alpha_n\eta_n^2\alpha_n^*\eta_n)\big)\leq \frac\pi2\tr\big(K\eta_n\gamma_n\eta_n\big).
\end{multline*}
Finally, the last term tends to zero by \eqref{estim_density_conc_compt} and thanks to the estimate $\|[\sqrt{K},\eta_n]\|\leq CR_n^{-1}$. In summary, we have derived the estimate (\ref{splitting}), as desired. 

Next, we notice $\cE(\gamma^1_n,\alpha_n^1)\geq I(\tr(\gamma_n^1))$ and $\cE(\gamma^2_n,\alpha_n^2)\geq I(\tr(\gamma_n^2))$. On the other hand, by Lemma \ref{lem:IHF},
\[ I(\lambda)\leq I(\tr(\gamma_n^1))+ I(\tr(\gamma_n^2)) +o(1), \]
using that $\tr(\gamma_n) \to \lambda$ and $\tr(\gamma_n)=\tr(\gamma^1_n)+\tr(\gamma^2_n)$, which follows from $\chi_n^2 + \zeta_n^2 \equiv 1$ and the cyclicity of the trace. Because of $\lim_{n\to\ii}\cE(\gamma_n,\alpha_n)=I(\lambda)$, we deduce
\begin{equation}
 \lim_{n\to\ii}\cE(\gamma^1_n,\alpha_n^1)= I(\lambda^1)\quad\text{and}\quad\lim_{n\to\ii}\cE(\gamma^2_n,\alpha_n)= I(\lambda-\lambda^1),
\label{CV_decomp}
\end{equation}
where we use that $\tr(\gamma^1_n)=\int\chi_n^2\rho_{\gamma_n}\to\lambda^1$, by Lemma \ref{lem:local_compactness}, and the the continuity of $\lambda\mapsto I(\lambda)$.

Therefore equation \eqref{CV_decomp} shows that both $\{ (\gamma_n^1,\alpha_n^1)\}$ and $\{(\gamma_n^2,,\alpha_n^2)\}$ are minimizing sequences for $I(\lambda^1)$ and $I(\lambda-\lambda^1)$, respectively. Since we already know that $\tr(\gamma^1_n)\to\lambda^1$, we conclude that $(\gamma_n^1,\alpha_n^1)\to(\gamma^1,\alpha^1)$ in $\cX$ from Corollary \ref{lem:compact}. By continuity of $\cE$, we also have $\lim_{n\to\ii}\cE(\gamma^1_n,\alpha_n^1)=\cE(\gamma^1,\alpha^1)=I(\lambda^1)$, so that $(\gamma^1,\alpha^1)$ is indeed a minimizer for $I(\lambda^1)$. This completes our proof of (\ref{eq:split}).

\subsubsection*{\bf Step 4: Binding Inequality and Conclusion}
At this stage, we have proven that the energy of our minimizing sequence behaves like
$$I(\lambda)= I(\lambda^1)+I(\lambda-\lambda^1)$$
and we have seen that there exists a minimizer $(\gamma^1,\alpha^1)$ for $I(\lambda^1)$. Note that $(\gamma^1,\alpha^1)$ is the weak limit of $(\gamma_n,\alpha_n)$ and the strong limit of $(\gamma_n^1,\alpha_n^1)=\chi_n(\gamma_n,\alpha_n)\chi_n$ in $\cX$. On the other hand, the sequence $\{ (\gamma_n^2,\alpha_n^2)\}$ is a minimizing sequence for $I(\lambda-\lambda^1)$. Hence one of the two following situations must occur:
\begin{itemize}
 \item \textbf{either} $\{ (\gamma_n^2,\alpha_n^2)\}$ is relatively compact up to translations. In this case it converges (up to a subsequence) to a minimizer $ (\gamma^2,\alpha^2)$ for $I(\lambda-\lambda^1)$;
\item \textbf{or} $\{ (\gamma_n^2,\alpha_n^2)\}$ is not relatively compact up to translations. We may then apply the whole process again and we deduce that we have
$$I(\lambda-\lambda^1)=I(\lambda^2)+I(\lambda-\lambda^1-\lambda^2),$$
for some $\lambda^2>0$ and that there exists a minimizer $(\gamma^2,\alpha^2)\in\cK_{\lambda^2}$ for $I(\lambda^2)$.
\end{itemize}
We can summarize both cases by saying that
\begin{equation} \label{eq:bindeq}
I(\lambda)=I(\lambda^1)+I(\lambda^2)+I(\lambda-\lambda^1-\lambda^2) ,
\end{equation}
for some $\lambda^1,\lambda^2>0$ such that $\lambda^1+\lambda^2\leq\lambda$ and that $I(\lambda^k)$ possess a minimizer $(\gamma^k,\alpha^k)$ for $k=1,2$.
Next, we claim the following result.
\begin{lemma}{\bf (Binding Inequality).} \label{lem:binding} Assume that $\lambda^1$ and $\lambda^2$ are as above. Then one has
\begin{equation}
 I(\lambda^1+\lambda^2)<I(\lambda^1)+I(\lambda^2).
\label{binding}
\end{equation}
\end{lemma}

\noindent Using \eqref{binding} and Lemma \ref{lem:IHF}, we deduce the strict inequality
\begin{equation}
I(\lambda)\leq I(\lambda^1+\lambda^2)+I(\lambda-\lambda^1-\lambda^2)<I(\lambda^1)+I(\lambda^2)+I(\lambda-\lambda^1-\lambda^2),
\end{equation}
which contradicts equation (\ref{eq:bindeq}). Hence the proof of Theorem \ref{thm:relaxed} will be finished after we have proven Lemma \ref{lem:binding}.


\begin{proof}[Proof of Lemma \ref{lem:binding}]
Consider two minimizers $(\gamma^1,\alpha^1)$ and $(\gamma^2,\alpha^2)$ for $I(\lambda^1)$ and $I(\lambda^2)$, respectively. We claim that   
\begin{equation} \label{eq:bind}
\cE(\chi_R\gamma^i \chi_R,\chi_R\alpha^i\chi_R)\leq I(\lambda^i)+ \frac{C}{R^2}, \quad \mbox{for $i=1,2$},
\end{equation}
for all $R>0$ sufficiently large and some constant $C > 0$, where we recall the definition of the cutoff functions $\chi_R$ and $\zeta_R$ from (\ref{eq:cutoff}).

It suffices to show (\ref{eq:bind}) for $i=1$, since the proof for $i=2$ is analogous. As for the kinetic energy, the localization formula in Lemma \ref{lem:IMS1} yields $\tr(T \gamma^1) \geq \tr(T \chi_R \gamma^1 \chi_R) - C/R^2$.
Next, we estimate the direct term as follows (using the symmetry in $x$ and $y$): 
\begin{multline*}
\big | D(\rho_{\gamma^1},\rho_{\gamma^1}) - D(\rho_{\chi_R \gamma^1 \chi_R}, \rho_{\chi_R \gamma^1 \chi_R}) \big | \\ \leq \iint_{\R^3 \times \R^3} \frac{(1-\chi_R^2(x) \chi_R^2(y)) \rho_{\gamma^1}(x) \rho_{\gamma^1}(y)}{|x-y|} \, dx \,dy \\
\leq 2 \iint_{\{ |x| \geq R \} \times \R^3} \frac{\rho_{\gamma^1}(x) \rho_{\gamma^1}(y)}{|x-y|} \, dx \, dy  \leq C \int_{|x| \geq R} \rho_{\gamma^1}(x) \, dx \leq \frac{C}{R^2},
\end{multline*}
for $R >0$ sufficiently large. Here we used the Kato-Hardy inequality to conclude that $\| |x|^{-1} \ast \rho_{\gamma^1} \|_{L^\infty} \leq C \|\sqrt{\rho_{\gamma^1}}\|_{H^{1/2}}$ as well as the decay estimate from Lemma \ref{lem:rho_decay}. 
Finally, we estimate the pairing term as follows: 
\begin{multline*}
 \iint_{\R^3 \times \R^3}  \frac{(1-\chi_R(x)^2 \chi_R(y)^2) |\alpha^1(x,y)|^2}{|x-y|} \, dx\, dy   \\ \leq 2 \iint_{\{ |x| \geq R\} \times \R^3} \frac{|\alpha^1(x,y)|^2}{|x-y|} \,dx \,dy  
   \leq 2 \iint_{\R^3 \times \R^3} \frac{\zeta_{R/2}(x)^2 |\alpha^1(x,y)|^2}{|x-y|} \, dx \, dy  \\ \leq  \pi  \tr (K \zeta_{R/2} \gamma^1\zeta_{R/2}) 
   \leq C \left (\int_{|x| \geq R/2} \rho_{\gamma^1} (x) \, dx + \tr(T \zeta_{R/2} \gamma^{1} \zeta_{R/2} ) \right ) .
\end{multline*}
Here the first inequality follows from the symmetry $|\alpha^1(x,y)|=|\alpha^1(y,x)|$ and the support properties of $\chi_R$. In the last line, we used Kato's inequality and $(\alpha^1)^* \alpha^1 \leq \gamma^1$. Putting all estimates together and using Lemma \ref{lem:rho_decay} again, we conclude that \eqref{eq:bind} holds. 

Next, we prove \eqref{binding}. To this end, we consider the trial state 
\begin{equation}
(\gamma_R,\alpha_R):=U(\chi_R\gamma^1\chi_R,\chi_R\alpha^1\chi_R)U^*+\tau_{5Rv}V(\chi_R\gamma^2\chi_R,\chi_R\alpha^2\chi_R)V^*\tau^*_{5Rv},
\end{equation}
for some unit vector $v$ and some rotations $U,V\in SO(3)$. As $\chi_R$ and $\tau_{5Rv}\chi_R$ have disjoint supports, it is easily seen that the above operator belongs to $\cK$. Moreover, the number of particles of this trial state satisfies
\begin{equation} 
\tr(\gamma_R)=\tr(\chi_R\gamma^1\chi_R)+\tr(\chi_R\gamma^2\chi_R)\leq \lambda^1+\lambda^2.
\end{equation}
Note the invariance by rotation and translation $\cE(\tau_vU(\gamma,\alpha)U^*\tau_v^*)=\cE(\gamma,\alpha)$ for all $v\in\R^3$ and all $U\in SO(3)$. Since $\lambda\mapsto I(\lambda)$ is decreasing by Lemma \ref{lem:IHF}, we have
\begin{align*}
I(\lambda^1+\lambda^2)&\leq \int_{SO(3)}\int_{SO(3)}\cE(\gamma_R,\alpha_R)\,dU\,dV\\
 &= \cE(\chi_R\gamma^1\chi_R,\chi_R\alpha^1\chi_R)+\cE(\chi_R\gamma^2\chi_R,\chi_R\alpha^2\chi_R)\\
&\qquad\qquad\qquad\qquad\qquad\qquad-\kappa\iint_{\R^3\times\R^3}\frac{f(x)g(y-5Rv)}{|x-y|}dx\,dy\\
 &\leq I(\lambda^1)+I(\lambda^2)+\frac{C}{R^2} - \frac{\kappa\lambda^1\lambda^2}{R},
\end{align*}
which follows from Newton's theorem and \eqref{eq:bind}, where we have defined
$$f(x)=\int_{SO(3)}\chi_R(x)^2\rho_{\gamma^1}(U\cdot x)\;dU,\quad g(y)=\int_{SO(3)}\chi_R(y)^2\rho_{\gamma^2}(V\cdot y)\;dV.$$
This last estimate yields the desired result by taking $R > 0$ sufficiently large.
\end{proof}

The proof of Theorem \ref{thm:relaxed} is now complete. \hfill $\blacksquare$

\section{Proof of Theorem \ref{thm:orbitals}}\label{sec:proof_thm_orbitals}
The proof of Theorem \ref{thm:orbitals} is a combination of ideas from the proof of Theorem \ref{thm:relaxed} and of the usual geometrical methods for $N$-body systems as was developed in \cite{Hun-66,VanWinter-64,Zhislin-60,Enss-77,Simon-77,Sigal-82}. To our knowledge, they were applied for the first time to nonlinear models approximating the $N$-body case (including the Hartree-Fock approximation) by Friesecke in \cite{Friesecke-03}. We shall only sketch the proof of Theorem  \ref{thm:orbitals}. In our opinion, it is much easier than the proof of Theorem \ref{thm:relaxed}.

\subsubsection*{\bf Step 1: First properties of HF energy} The first step consists in showing that
\begin{equation}
 \IHF(1)=0,\qquad  \IHF(N)<0\text{ for all $N\geq2$}
\label{prop_IHF_1}
\end{equation}
and that
\begin{equation}
 \IHF(N)\leq \IHF(N-K)+\IHF(K)\text{ for all } K=1..N-1.
\label{prop_IHF_2}
\end{equation}
The value of $\IHF(1)$ is obvious as for density matrices of rank one, the exchange and direct terms cancel. The rest is proved similarly as in Lemma \ref{lem:IHF}.

\medskip

The main part of the proof consists in showing that if 
\begin{equation}
 \IHF(N)< \min(\IHF(N-K)+\IHF(K),\ K=1..N-1).
\label{binding_HF}
\end{equation}
holds true, then all minimizing sequences are compact and thus one gets the existence of a minimizer for $\IHF(N)$. In the following we fix some $N\geq2$, assume that \eqref{binding_HF} holds, and consider a minimizing sequence $\{\gamma_n\}$ for  $\IHF(N)$. Similarly to Lemma \ref{lem:coercive}, we obtain that $\{(\gamma_n,0)\}$ is bounded in $\cX$. 

\subsubsection*{\bf Step 2: No Vanishing} We claim that our sequence $\{\gamma_n\}$ cannot vanish, similarly to Lemma \ref{lem:novanish}. Indeed if  $\{\gamma_n\}$ would vanish, applying Lemma \ref{lem:vanish} we get that $D(\rho_{\gamma_n},\rho_{\gamma_n})\to0$, hence we would obtain $\IHF(N)\geq0$ which contradicts \eqref{prop_IHF_1}. Therefore we may assume (translating $\gamma_n$ if necessary), that $\gamma_n\wto\gamma$ with $\lambda^1:=\tr(\gamma)>0$. Also we have similarly to Lemma \ref{lem:local_compactness} that we can find a cut-off function on a ball of size $R_n$ such that $\sqrt{K}\chi(R_n^{-1}\cdot)\gamma_n\chi(R_n^{-1}\cdot)\sqrt{K}\to\sqrt{K}\gamma\sqrt{K}$ strongly.
Also, similarly to Corollary \ref{lem:compact}, one shows that $\{\gamma_n\}$ is compact when $\tr(\gamma)=N$, hence it suffices to prove that $\lambda^1=N$.

\subsubsection*{\bf Step 3: $N$-body Geometrical Methods} The next step is an easy adaptation of \cite{Friesecke-03}, i.e. it consists in using a cluster decomposition like in the usual HVZ Theorem, based on the fact that $\gamma_n$ arises from an $N$-body wavefunction which is a Slater determinant. Indeed following \cite{Friesecke-03}, one arrives at an estimate of the form (up to extraction of a subsequence):
\begin{equation}
\sum_{K=1}^{N-1}\bigg\{(\IHF(N)-\IHF(N-K)-\IHF(K)\bigg\}\lim_{n\to\ii}z_K^n\geq0.
\label{esztim_z_n}
\end{equation}
The nonnegative numbers $z_K^n$ are defined by
$$z_K^n=\left(\begin{matrix}N\\ K\end{matrix}\right)\int\cdots\int \prod_{i=1}^K\zeta_n(x_i)^2\prod_{i=K+1}^N\chi_n(x_i)^2|\Psi_n(x_1,...,x_N)|^2dx_1\cdots dx_N,$$
where $\Psi_n$ is the Slater determinant associated with $\gamma_n$ and $\chi_n=\chi(R_n^{-1}\cdot)$ as introduced before, $\zeta_n=\sqrt{1-\chi_n^2}$, see \cite{Friesecke-03}. Note that $z_K^n$ may be interpreted as the norm of the wavefunction associated with all the clusters for which $K$ particles are sent to infinity while $N-K$ particles stay close to zero.
Note one has
\begin{equation}
\sum_{K=0}^{N}z_K^n=1\ \text{ and }\ \lim_{n\to\ii}\sum_{K=0}^{N}(N-K)z_K^n=\lambda^1.
\label{relations_z_n} 
\end{equation}
The strict inequality \eqref{binding_HF} and \eqref{esztim_z_n} imply that we must have 
$$\lim_{n\to\ii}z_K^n=0\ \text{ for } K=1..N-1.$$
This implies that we have
$$\lim_{n\to\ii}\int\cdots\int \zeta_n(x_1)^2\chi_n(x_2)^2|\Psi_n(x_1,...,x_N)|^2dx_1\cdots dx_N=0,$$
$$\lim_{n\to\ii}\norm{\Psi_n-\sqrt{z_0^n}\Psi^1_n-\sqrt{z_N^n}\Psi^2_n}=0,$$
$$\lim_{n\to\ii}\norm{\gamma_{\Psi_n}-z_0^n\gamma_{\Psi_n^1}-z_N^n\gamma_{\Psi_n^2}}_{\gS_1}=\lim_{n\to\ii}\norm{\gamma_{\Psi_n^1}\gamma_{\Psi_n^2}}_{\gS_1}=0$$
where $\Psi_n^1$ and $\Psi_n^2$ are Hartree-Fock states defined as
$$\Psi_n^1:=\frac{\chi_n^{\otimes N}\Psi_n}{\norm{\chi_n^{\otimes N}\Psi_n}},\qquad\Psi_n^2:=\frac{\zeta_n^{\otimes N}\Psi_n}{\norm{\zeta_n^{\otimes N}\Psi_n}}.$$
Using that the operators $\gamma_{\Psi_n}$, $\gamma_{\Psi_n^1}$ and $\gamma_{\Psi_n^2}$ are projectors, we obtain
$$z_0^n\gamma_{\Psi_n^1}+z_N^n\gamma_{\Psi_n^2}+o(1)=\gamma_{\Psi_n}=\big(\gamma_{\Psi_n}\big)^2=(z_0^n)^2\gamma_{\Psi_n^1}+(z_N^n)^2\gamma_{\Psi_n^2}+o(1).$$
Taking the trace of the previous equality, we infer
$$\lim_{n\to\ii}\big(z_0^n-(z_0^n)^2\big)=\lim_{n\to\ii}\big(z_N^n-(z_N^n)^2\big)=0.$$
By \eqref{relations_z_n}, we have $\lim_{n\to\ii}Nz_0^n=\lambda^1>0$, therefore we arrive at the following conclusion: 
$$\lim_{n\to\ii}z_0^n=1\ \text{ and }\ \lim_{n\to\ii}z_K^n=0\ \text{ for } K=1..N$$
hence $\lambda^1=N$. This proves that $\{(\gamma_n,0)\}$ is indeed compact in $\cX$.

\subsubsection*{\bf Step 4: Proof of the binding inequality \eqref{binding_HF}} 
We have seen that minimizing sequences are all compact when \eqref{binding_HF} hold. It remains to prove \eqref{binding_HF}. This is done like in \cite{Friesecke-03} by induction on $N\geq2$. As $\IHF(1)=0$ and $\IHF(2)<0=2\IHF(1)$, we deduce that $\IHF(2)$ possesses a minimizer. Next the strict inequality $\IHF(3)<\IHF(2)+\IHF(1)=\IHF(2)$ is proved as in \cite[Thm. 3]{Friesecke-03}, using that $\IHF(2)$ possesses a minimizer. The rest follows by induction, as soon of course as the energy stays bounded from below and coercive, i.\,e.~when $N<N^{\rm HF}(\kappa)$ holds. 

\medskip
This completes our sketch of the proof of Theorem \ref{thm:orbitals}. \hfill$\blacksquare$\\

\begin{appendix}

\section{Localization of Kinetic Energy}\label{appendix:localization}

In the section we collect some localization estimates for the pseudo-relativistic kinetic energy operator $\sqrt{-\Delta + m^2}$. In fact, similar localization estimates -- often referred to as {\em IMS-type localization formulas} --  can be found, for instance, in the stability analysis of relativistic matter \cite{LieYau-88}. However, the following IMS-type estimates for $\sqrt{-\Delta + m^2}$ turn out to be suitable when studying lack of compactness of minimizing sequences. 

\begin{lemma} \label{lem:IMS1}
Let $\{ \chi_i \}_{i \in I}$ be a smooth partition of unity with $\sum_{i \in I} \chi_i(x)^2 \equiv 1$ such that each $\nabla \chi_i \in L^p(\R^d)$ with $p \in (2d, \infty]$. Then
\begin{equation}
K \geq \sum_{i \in I} \chi_i K \chi_i - \frac{1}{\pi} \int_0^\infty \frac{1}{s + K^2} \left ( \sum_{i \in I} |\nabla \chi_i |^2 \right ) \frac{1}{s + K^2} \, \sqrt{s} \, ds,
\label{eq:localization_formula}
\end{equation} 
where $K = \sqrt{-\Delta + m^2}$ with $m >0$.
\end{lemma}

The proof of Lemma \ref{lem:IMS1} readily follows from the next result.

\begin{lemma} \label{lem:IMS2}
Suppose $\chi \in C^\infty(\R^d)$ such that $\nabla \chi \in L^p(\R^d)$ with $p \in (2d,\infty]$, and let $K = \sqrt{-\Delta + m^2}$ with $m >0$. Then the following formula holds:
\begin{equation} \label{eq:lemloc}
\frac{1}{2} ( K \chi^2 + \chi^2 K ) = \chi K \chi - \frac{1}{\pi} \int_0^\infty \frac{1}{s+K^2} | \nabla \chi |^2 \frac{1}{s+K^2} \, \sqrt{s} \, ds + L_\chi,
\end{equation}
where $L_\chi$ is a nonnegative bounded operator. Moreover, we have
\begin{equation}
\Big \| \int_0^\infty \frac{1}{s+K^2} |\nabla \chi|^2 \frac{1}{s+K^2} \, \sqrt{s} \, ds \Big \|_{\gS_{p/2}(L^2(\R^d))} \leq C \| \nabla \chi \|_{L^p(\R^d)}^2,
\end{equation}
\begin{equation}
\| L_\chi \|_{\gS_{p/2}(L^2(\R^d))} \leq C \| \nabla \chi \|_{L^p(\R^d)}^2,
\label{estim_L_localization}
\end{equation}
for some constant $C$ independent of $\chi$ (but depending on $m$, $p$ and $d$).
\end{lemma}

\begin{proof}
By functional calculus, we have the representation formula
\begin{equation} \label{eq:Krep}
K = \frac{1}{\pi} \int_0^\infty \frac{K^2}{s+K^2} \frac{ds}{\sqrt{s}} ,
\end{equation}
which leads to
\begin{equation} \label{eq:Kcom}
[K, B] = \frac{1}{\pi} \int_0^\infty \frac{1}{s+K^2} [K^2, B ] \frac{1}{s+K^2} \sqrt{s} \, ds .
\end{equation}
Using (\ref{eq:Kcom}), a calculation shows
\begin{align*}
\frac{1}{2} ( K \chi^2 + \chi^2 K ) & = \chi K \chi + \frac{1}{2} [[K,\chi],\chi] \\
& = \chi K \chi + \frac{1}{2 \pi} \int_0^\infty\left [  \frac{1}{s+K^2} [ K^2, \chi ] \frac{1}{s+K^2}, \chi \right ] \, \sqrt{s} \, ds \\
& = \chi K \chi + \frac{1}{2 \pi} \int_0^\infty \frac{1}{s+K^2} \big [ [K^2, \chi ], \chi \big ] \frac{1}{s+K^2}  \sqrt{s} \, ds + L_\chi, 
\end{align*}
with
\begin{equation}
L_\chi = \frac{1}{\pi} \int_0^\infty \frac{1}{s+K^2} [K^2,\chi] \frac{1}{s+K^2} [\chi, K^2] \frac{1}{s+K^2} \, \sqrt{s} \, ds .
\end{equation}
Since $[[K^2,\chi],\chi] = - 2 |\nabla \chi|^2$, this yields (\ref{eq:lemloc}). Furthermore, we note that
\begin{equation}
L_\chi = \frac{1}{\pi} \int_0^\infty A(s)^* A(s) \sqrt{s} \, ds, \quad A(s) = \frac{1}{(s+K^2)^{1/2}} [\chi, K^2] \frac{1}{s+K^2},
\end{equation}
which shows that $L_\chi$ is nonnegative.

We now show that all the above commutators are bounded (which will also validate our somewhat formal calculation). Recall the Kato-Seiler-Simon inequality (see \cite{SeiSim-75} and \cite[Thm 4.1]{Simon-79})
\begin{equation}
 \norm{f(-i \nabla_x)g(x)}_{\gS_p(L^2(\R^d))}\leq (2\pi)^{-d/p}
\norm{f}_{L^p(\R^d)} \norm{g}_{L^p(\R^d)} , \quad \mbox{for $p \in [2,\infty]$}.
\label{KSS}
\end{equation}
This yields, choosing some $d/(2p)<q<1$, the following bound. For all $s \geq 0$,
\begin{multline}
\norm{\frac{1}{s+K^2}|\nabla\chi|^2\frac{1}{s+K^2}}_{\gS_{p/2}(L^2(\R^d))}\\
\leq \frac{C}{(m^2+s)^{2-2q}}\norm{(s+\xi^2)^{-q}}_{L^p(\R^d)}^2\norm{\nabla\chi}_{L^p_\xi(\R^d)}^2.
\end{multline}
Furthermore, a scaling argument shows that 
$$\norm{(s+\xi^2)^{-q}}_{L^p_\xi(\R^d)}=s^{d/(2p)-q}\norm{(1+\xi^2)^{-q}}_{L^p_\xi(\R^d)} .$$
Hence we obtain
\begin{multline}
\norm{\int_0^\ii\frac{1}{s+K^2}|\nabla\chi|^2\frac{1}{s+K^2}\sqrt{s}\,ds}_{\gS_{p/2}(L^2(\R^d))}\\
\leq C\norm{\nabla\chi}_{L^p(\R^d)}^2\int_0^\ii\frac{s^{1/2+d/p-2q}}{(m^2+s)^{2-2q}}\,ds ,
\end{multline}
where the last integral is convergent under our assumptions on $d$ and $p$.
The argument is exactly the same for $L_\chi$, noting that $[K^2,\chi]=-i(P\cdot\nabla\chi+\nabla\chi\cdot P)$ with $P= -i \nabla$ and using that $\norm{P(s+m^2+|P|^2)^{-1/2}}\leq 1$. \end{proof}

\section{Proof of Proposition \ref{prop:def_max_mass}} \label{appendix:proof_def_max_mass} 
We consider the following minimization problem
\begin{multline}
I_{0,\kappa}(\lambda)=\inf_{(\gamma,\alpha)\in\cK_\lambda}\left\{\tr(\sqrt{-\Delta}\gamma) -\frac{\kappa}{2}D(\rho_\gamma,\rho_\gamma) \right . \\ \left  . +\frac\kappa2 \Ex(\gamma)-\frac\kappa2\iint\frac{|\alpha(x,y)|^2}{|x-y|}dx\,dy\right\}.
\label{def:I_0}
\end{multline}
By rescaling $\gamma(x,y) \mapsto \beta^{3} \gamma(\beta x, \beta y)$ and $\alpha(x,y) \mapsto \beta^{3} \alpha(\beta x, \beta y)$ with $\beta > 0$, we deduce that either $I_{0,\kappa}(\lambda)=0$ or $-\ii$. Following the proof of Lemma \ref{lem:IHF}, one also sees that $I_{0,\kappa}$ is non-increasing. Hence there is a uniquely defined $\lambda^{\rm HFB}(\kappa)$ such that $I_{0,\kappa}(\lambda)=0$ if and only if $\lambda\leq \lambda^{\rm HFB}(\kappa)$. 
Next, we notice that
$$\sqrt{-\Delta}-m\leq \sqrt{-\Delta+m^2}-m\leq \sqrt{-\Delta}$$
hence
$I_{0,\kappa}(\lambda)-m\lambda\leq I(\lambda)\leq I_{0,\kappa}(\lambda)$.
This shows that $ I(\lambda)$ is bounded below if and only if $\lambda\leq \lambda^{\rm HFB}(\kappa)$. Furthermore, it is clear that $\lambda^{\rm HFB}(\cdot)$ is nonincreasing, since for all $(\gamma,\alpha)\in\cK$,
$$-D(\rho_\gamma,\rho_\gamma)+ \Ex(\gamma)-\iint\frac{|\alpha(x,y)|^2}{|x-y|}dx\,dy\leq0$$
by \eqref{estim_exchange}. The continuity of $\lambda^{\rm HFB}(\kappa)$ is easy to verify. 

Next, we prove that $\lambda^{\rm HFB}(\kappa)>0$ and we derive the asymptotic behavior of $\lambda^{\rm HFB}(\kappa)$ as $\kappa\to 0$ stated in Proposition \ref{prop:def_max_mass}. To this end, we introduce another homogeneous minimization problem
\begin{equation}
I_{0,\kappa}^{\rm red}(\lambda)=\inf_{(\gamma,0)\in\cK_\lambda}\left\{\tr(\sqrt{-\Delta}\gamma) -\frac{\kappa}{2}D(\rho_\gamma,\rho_\gamma)\right\}.
\end{equation}
which, of course, yields a corresponding $\lambda^{\rm red}(\kappa)$. We claim that $\lambda^{\rm HFB}(\kappa)\sim_{\kappa\to0}\lambda^{\rm red}(\kappa)$. Indeed, by the nonnegativity of the exchange term and \eqref{estim_pairing_term},
\begin{equation}
I_{0,\kappa}(\lambda)\geq \left(1-\frac{\kappa\pi}{4}\right)I_{0,\frac{\kappa}{1-\kappa\pi/4}}^{\rm red}(\lambda).
\end{equation}
Hence
\begin{equation}
\lambda^{\rm HFB}(\kappa)\geq \lambda^{\rm red}\left(\frac\kappa{1-\kappa\pi/4}\right).
\label{estim_below_reduced}
\end{equation}

Now we recall that $\lambda^{\rm red}(\kappa')>0$ for every $\kappa'>0$. Indeed we have 
\begin{equation}
\iint_{\R^3 \times \R^3} \frac{\rho_\gamma(x) \rho_\gamma(y)}{|x-y|} dx \, dy \leq C' \| \rho_\gamma \|_{L^1}^{2/3} \| \rho_\gamma \|_{L^{4/3}}^{4/3}\leq C\lambda^{2/3}\tr (\sqrt{-\Delta} \gamma)
\end{equation}
where we have used the Hardy-Littlewood-Sobolev inequality and the well-known \cite{Daubechies-83} semi-classical estimate $\tr (\sqrt{-\Delta} \gamma) \geq K \| \rho_\gamma \|_{L^{4/3}}^{4/3}$ for fermionic density matrices $0\leq\gamma\leq1$. This shows that 
$\lambda^{\rm red}(\kappa')\geq \big(2/(C\kappa')\big)^{3/2}>0$. Hence we deduce from \eqref{estim_below_reduced} and $0\leq\kappa<4/\pi$ that $\lambda^{\rm HFB}(\kappa)>0$.

Next we derive an upper bound for $\lambda^{\rm HFB}(\kappa)$. We have for any $(\gamma,0)\in\cK_\lambda$,
\begin{equation}
I_{0,\kappa}(\lambda)\leq \cE(\gamma,0)=\tr(\sqrt{-\Delta}\gamma) -\frac{\kappa}{2}D(\rho_\gamma,\rho_\gamma)+\frac\kappa2 \Ex(\gamma).
\end{equation}
Now we have similarly to \eqref{estim_pairing_term}
\begin{equation}
\iint\frac{|\gamma(x,y)|^2}{|x-y|}dx\,dy\leq \frac\pi2\tr(\sqrt{-\Delta}\gamma^2)\leq \frac\pi2\tr(\sqrt{-\Delta}\gamma) .
\end{equation}
Hence we obtain
\begin{equation}
I_{0,\kappa}(\lambda)\leq \left(1+\frac{\kappa\pi}{4}\right)I_{0,\frac{\kappa}{1+\kappa\pi/4}}^{\rm red}(\lambda),
\end{equation}
and we thus get the upper bound
\begin{equation}
\lambda^{\rm HFB}(\kappa)\leq \lambda^{\rm red}\left(\frac\kappa{1+\kappa\pi/4}\right).
\label{estim_above_reduced}
\end{equation}

The last step consists in showing that $\lambda^{\rm red}(\kappa)$ behaves as stated when $\kappa \to 0$, which will imply the result for $\lambda^{\rm HFB}(\kappa)$ by \eqref{estim_below_reduced} and \eqref{estim_above_reduced}. This was indeed essentially done by Lieb and Yau, see Corollary 1 in \cite{LieYau-87}.
\hfill$\blacksquare$\\

\section{Proof of Proposition \ref{prop:comparison_I_J}}\label{appendix:F}
We want to prove that $G(\lambda)=\beta\lambda$, where $\beta$ is the first eigenvalue of the operator $K-\kappa/2|x|$ acting on $L^2_{\rm odd}(\R^3)$ when $q=1$ and $L^2(\R^3)$ when $q\geq2$. For the reader's convenience, we recall that $K =\sqrt{-\Delta + m^2}$.

We start by deriving a lower bound. Using \eqref{estim_pairing}, we infer $\tr(K\gamma)\geq\tr(K\alpha\alpha^*)$ and thus
\begin{equation}
G(\lambda)\geq \inf_{\substack{\alpha^T=-\alpha,\\ \tr(\alpha^*\alpha)\leq\lambda}}\left\{\tr(K\alpha\alpha^*)-\frac{\kappa}{2}\iint_{\R^3\times\R^3}\frac{|\alpha(x,y)|^2}{|x-y|}dx\,dy\right\} . 
\end{equation}
The right-hand side of the previous inequality can also be written, using the antisymmetry of $\alpha$, as follows.
\begin{equation}
\inf_{\substack{\alpha(x,y)^T=-\alpha(x,y),\\ \norm{\alpha}_{L^2(\R^3 \times \R^3)}^2\leq \lambda}}\pscal{\left(\frac{K_x+K_y}{2}-\frac{\kappa}{2|x-y|}\right)\alpha,\alpha} .
\end{equation}
This quantity is just $\lambda$ times the lowest eigenvalue of the two-body operator
\begin{equation}
\frac{K_x+K_y}{2}-\frac{\kappa}{2|x-y|}
\end{equation}
when it is restricted to the fermionic space $L^2(\R^3,\C^q)\wedge L^2(\R^3,\C^q)$. 
If $q=1$, its lowest eigenvalue equals the first eigenvalue of the same operator acting on $L^2(\R^3_x,\R)\wedge L^2(\R^3_y,\R)$ (triplet state), whereas when $q\geq2$, it coincides with the first eigenvalue of the same operator acting on $L^2(\R^3,\R)\otimes_s L^2(\R^3,\R)$ (singlet state).
Removing the center of mass as explained in \cite{LewSieVug-97}, one sees that this lowest eigenvalue is given by
$$\inf_{\ell\in\R^3}\inf\sigma_{\gH}\left(\frac12\sqrt{\left|-i\nabla_{u}+\frac{\ell}{2}\right|^2+m^2}+\frac12\sqrt{\left|-i\nabla_{u}-\frac{\ell}{2}\right|^2+m^2}-\frac{\kappa}{2|u|}\right), $$
where $u=x-y$ and $\ell$ is the Fourier variable associated with the space variable $v=(x+y)/2$. We have to take $\gH=L^2_{\rm odd}(\R^3)$ when $q=1$ and $\gH=L^2(\R^3)$ when $q\geq2$. The infimum above is attained for $\ell=0$, and hence we obtain the desired lower bound:
\begin{equation}
G(\lambda)\geq \lambda\inf\sigma_{\gH}\left(\sqrt{-\Delta_u+m^2}-\frac{\kappa}{2|u|}\right):=\beta\lambda.
\end{equation}
Note that we can in fact deduce that the infimum on the right side is attained. This follows from the min-max principle and the operator inequality $\sqrt{-\Delta + m^2} - \frac{\kappa}{2 |x|} \leq -\frac{1}{2m}\Delta + m - \frac{\kappa}{2|x|}$, where the right side is the Schr\"odinger operator for the nonrelativistic hyrdogen atom.

Next, we prove the upper bound in a similar way as done for the M\"uller functional in \cite{FraLieSeiSie-07}. Since $K-\kappa/2|x|$ is a real operator, there is a real-valued eigenfunction $f \in \gH$ such that 
\begin{equation}
\left(\sqrt{-\Delta + m^2} - \frac{\kappa}{2 |x|}\right) f = \beta f .
\end{equation}
Using the method in \cite{DalOstSto-08}, one can show that $f\in C^\ii(\R^3\setminus\{0\})$ and that it decays exponentially at infinity when $0\leq\kappa<4/\pi$. Note that $f\in D\big(\sqrt{-\Delta + m^2} - \kappa/(2 |x|)\big)\subset H^{1/2}(\R^3)$ but this domain is only known to be $H^1(\R^3)$ when $\kappa\leq1$. (We will address this regularity issue below again).

Next, we fix a nonnegative function $\chi \in C^\infty_0(\R^3)$  such that $\int\chi^4=\lambda$, and we denote $\chi_L(x)=L^{-3/4}\chi(x/L)$ for $L >0$ given. We define
\begin{equation}
\alpha_L(x,y)=\chi_L(x)f(x-y)\chi_L(y)\Sigma,
\end{equation}
where $\Sigma=1$ if $q=1$ and, when $q\geq2$,
$$\Sigma=\frac{1}{\sqrt{2}}\left(\begin{matrix}
0&1&&&\\
-1&0&&&\\
&&0&&\\
&&&\ddots&\\
&&&&0
\end{matrix}\right).
$$
As $\alpha_L\in L^2(\R^3\times\R^3,\C^q\times\C^q)$, we may denote by $\alpha_L$ the Hilbert-Schmidt operator whose kernel is $\alpha_L(x,y)$. Then $\alpha_L^T=-\alpha_L$ by the choice of $f$. Note that
\begin{equation}
\norm{\alpha_L}\leq C\norm{f}_{L^1}\norm{\chi_L}_{L^\ii}^2=O(L^{-3/2}),
\end{equation}
and hence $\alpha_L\to0$ in the operator norm. Next, we define $\gamma_L$ as the unique nonnegative trace-class operator solving the following equation for $L>0$ large enough:
\begin{equation}
\gamma_L(1-\gamma_L)=\alpha_L\alpha_L^*. 
\label{equation_gamma_alpha}
\end{equation}
Note that $\gamma_L\leq 2\alpha_L\alpha_L^*$ if $L > 0$ is sufficiently large. In particular, we have $\norm{\gamma_L}=O(L^{-3})$. In analogy to \cite{FraLieSeiSie-07}, it can be checked using $\int\chi^4=\lambda$ that
$\tr(\alpha_L\alpha_L^*)=\lambda+O(L^{-2})$.
Therefore, we have $\tr(\gamma_L)=\lambda+O(L^{-2})$, by \eqref{equation_gamma_alpha}, and the fact that $\tr(\gamma_L^2)=O(L^{-3})$.

Next, we show that $(\gamma_L,\alpha_L)\in\cK$. To see this, we recall from \cite{Coleman-63} and \cite[Prop. 2]{Lewin-04a} that there exists an orthonormal basis $(\phi_i)$ of $L^2(\R^3,\C^q)$ such that
\begin{equation} \label{eq:feigen}
\alpha_L=\sum_{i\geq1}c_i\; \phi_{2i-1}\wedge\phi_{2i}.
\end{equation}
On the other hand, the operator $\alpha_L\alpha_L^*$ is just $1/2$ times the one-body density matrix of the two-body wavefunction $\alpha_L$, so that
\begin{equation}
\alpha_L\alpha_L^*=\frac12\sum_{i\geq1}(c_i)^2\left(|\phi_{2i-1}\rangle\langle\phi_{2i-1}|+|\phi_{2i}\rangle\langle\phi_{2i}|\right ),
\end{equation}
\begin{equation}
\gamma_L=\sum_{i\geq1}\gamma_i\left(|\phi_{2i-1}\rangle\langle\phi_{2i-1}|+|\phi_{2i}\rangle\langle\phi_{2i}|\right) ,
\end{equation}
where for all $i\geq1$, $\gamma_i(1-\gamma_i)=(c_i)^2/2$. Next, we have to verify that 
\begin{equation}
0\leq\Gamma_L=\left(\begin{matrix}
\gamma_L&\alpha_L\\ \alpha_L^*&1-\gamma_L
\end{matrix}\right)\leq1.
\end{equation}
As $\alpha_L\phi_{2i-1}=-\frac{c_i}{\sqrt{2}}\phi_{2i}$ and $\alpha_L\phi_{2i}=\frac{c_i}{\sqrt{2}}\phi_{2i-1}$ with respect to the basis $(\phi_i)$, this amounts to studying the following matrix
\begin{equation}
M=\left(\begin{matrix}
\gamma_i & 0 & 0 & c_i/\sqrt{2}\\
0 & \gamma_i & -c_i/\sqrt{2} &0\\
 0 & -c_i/\sqrt{2}&1-\gamma_i & 0\\
c_i/\sqrt{2} & 0 &0 & 1-\gamma_i 
\end{matrix}\right) ,
\end{equation}
which is easily shown to satisfy $0\leq M\leq 1$ (it is indeed a projection). Thus we have proven that $(\gamma_L,\alpha_L)\in\cK$ for $L > 0$ sufficiently large.
Next, we estimate
\begin{equation}
\cG(\gamma_L,\alpha_L)=\tr(K\alpha_L\alpha_L^*)-\frac\kappa2\iint\frac{|\alpha_L(x,y)|^2}{|x-y|}dx\,dy+\tr(K\gamma_L^2) .
\end{equation}
We start by discussing the kinetic energy. Using \eqref{eq:lemloc}, we obtain
\begin{multline*}
 \tr(K\alpha_L\alpha_L^*) \leq (2\pi)^{3/2}\tr\left(\chi_L^2K(p)\widehat{f}(p)\chi_L^2\widehat{f}(p)\right)\\
+\frac{(2\pi)^{3/2}}{\pi} \int_0^\ii\sqrt{s}\,ds \tr\left( \frac{\widehat{f}(p)}{K(p)^2+s}|\nabla\chi_L|^2 \frac{\widehat{f}(p)}{K(p)^2+s}\chi_L^2\right)
\end{multline*}
Since the eigenfunction $f$ satisfies (\ref{eq:feigen}), we obtain
$$(2\pi)^{3/2}\tr\left(\chi_L^2K(p)\widehat{f}(p)\chi_L^2\widehat{f}(p)\right)=\frac\kappa2\iint\frac{|\alpha_L(x,y)|^2}{|x-y|}dx\,dy+\beta\tr(\alpha_L\alpha_L^*).$$
Notice the above expressions are well-defined even when $f$ only belongs to $H^{1/2}(\R^3)$. When $4/\pi > \kappa>1$, they can be justified by approximating the eigenfunction $f$ in $H^{1/2}(\R^3)$ by smooth functions and passing to the limit. On the other hand, we have
\begin{align*}
 \norm{\frac{\widehat{f}(p)}{K(p)^2+s}|\nabla\chi_L|^2 \frac{\widehat{f}(p)}{K(p)^2+s}\chi_L^2}_{\gS_1}&\leq \norm{\frac{\widehat{f}(p)}{K(p)^2+s}|\nabla\chi_L|^2}_{\gS_2}\norm{\frac{\widehat{f}(p)}{K(p)^2+s}\chi_L^2}_{\gS_2}\\
&\leq \frac{\norm{f}_{L^2}\norm{\nabla\chi_L}_{L^4}^2}{m+s}\times \frac{\norm{f}_{L^2}\norm{\chi_L}_{L^4}^2}{m+s}\leq \frac{CL^{-2}}{(m+s)^2}.
\end{align*}
Moreover, for $L >0$ sufficiently large, 
\begin{equation}
\tr(K\gamma_L^2)\leq 2\tr\left(K(\alpha_L\alpha_L^*)^2\right)\leq 2\norm{\alpha_L}^2\tr(K\alpha_L\alpha_L^*)=O(L^{-3}). 
\label{estim_kinetic}
\end{equation}
In summary, this shows $\cG(\gamma_L,\alpha_L)\leq \beta\lambda+O(L^{-2})$.
Passing to the limit $L\to\ii$, we get the upper bound $G(\lambda)\leq \beta\lambda$. Combined with the lower bound, this finally yields the desired equality $G(\lambda)=\beta\lambda$.

The estimate for $I(\lambda)$ is obtained using the same state $(\gamma_L,\alpha_L)$. We have
$$\iint_{\R^3\times\R^3}\frac{|\gamma_L(x,y)|^2}{|x-y|}dx\,dy\leq \frac\pi2\tr(K\gamma_L^2)=O(L^{-3}).$$
which we proved above in \eqref{estim_kinetic}. Moreover note that $\rho_{\gamma_L}\geq \rho_{\alpha_L\alpha_L^*}$ by \eqref{equation_gamma_alpha}, and this is seen to imply
$$D(\rho_{\gamma_L},\rho_{\gamma_L})\geq D(\rho_{\alpha_L\alpha^*_L},\rho_{\alpha_L\alpha^*_L})=\frac1L D(\chi^4,\chi^4)+o(L^{-1})$$
as can be deduced from \cite[p. 052517]{FraLieSeiSie-07}. Therefore, by choosing $L >0$ sufficiently large, we obtain
$$
I(\lambda)\leq G(\lambda)-m\lambda-\frac{\kappa}{2L}D(\chi^4,\chi^4)+o(L^{-1})< G(\lambda)-m\lambda,
$$
as desired.\hfill$\blacksquare$\\

\end{appendix}


\end{document}